# TIME-REVERSAL AND ELLIPTIC BOUNDARY VALUE PROBLEMS

By Zhen-Qing Chen[1] and Tusheng Zhang

*University of Washington and University of Manchester*

In this paper, we prove that there exists a unique, bounded continuous weak solution to the Dirichlet boundary value problem for a general class of second-order elliptic operators with singular coefficients, which does not necessarily have the maximum principle. Our method is probabilistic. The time reversal of symmetric Markov processes and the theory of Dirichlet forms play a crucial role in our approach.

**1. Introduction.** The pioneering work by Kakutani [18] in 1944 on representing the solution to the classical Dirichlet boundary value problem

$$\begin{cases} \Delta u = 0, & \text{in } D, \\ u = f, & \text{on } \partial D, \end{cases}$$

using Brownian motion started a new era in the very fruitful interplay between probability theory and analysis. Here $D$ is a bounded connected open subset of $\mathbb{R}^n$. Since then, in place of Laplacian $\Delta$, there are two classes of second-order elliptic differential operators that have been studied in connection with probabilistic approach. One is the nondivergence form operator

$$(1.1) \qquad \mathcal{L} = \frac{1}{2} \sum_{i,j=1}^n a_{ij}(x) \frac{\partial^2}{\partial x_i \partial x_j} + \sum_{i=1}^n b_i(x) \frac{\partial}{\partial x_i} + q.$$

The other is the divergence form operator

$$(1.2) \qquad \mathcal{L} = \frac{1}{2} \sum_{i,j=1}^n \frac{\partial}{\partial x_i} \left( a_{ij}(x) \frac{\partial}{\partial x_j} \right) + \sum_{i=1}^n b_i(x) \frac{\partial}{\partial x_i} + q,$$

Received March 2008; revised August 2008.

[1]Supported in part by NSF Grant DMS-06-00206.

*AMS 2000 subject classifications.* Primary 60J70, 60J57; secondary 35R05, 31C25i, 60H05, 60G46.

*Key words and phrases.* Diffusion, time-reversal, Girsanov transform, Feynman–Kac transform, multiplicative functional, partial differential equation, weak solution, boundary value problem, quadratic form, probabilistic representation.

where $A(x) = (a_{ij}(x))$ is an $n \times n$ symmetric bounded positive definite matrix. For nondivergence form operator $\mathcal{L}$ in (1.1), one can run stochastic differential equation

$$dX_t = \sigma(X_t)\, dB_t + \mathbf{b}(X_t)\, dt,$$

where $\sigma$ is a symmetric $n \times n$ matrix such that $\sigma^2 = A$, $\mathbf{b} = (b_1, \ldots, b_n)$ and $B$ is a Brownian motion on $\mathbb{R}^n$. The infinitesimal generator of $X$ is $\mathcal{L}_0 = \mathcal{L} - q$. Under some suitable conditions, the solution for $\mathcal{L}u = 0$ in $D$ with $u = f$ on $\partial D$ can be solved by

$$u(x) = \mathbf{E}_x\left[\exp\left(\int_0^{\tau_D} q(X_s)\, ds\right) f(X_{\tau_D})\right] \qquad \text{for } x \in D;$$

see [12]. Here $\tau_D = \inf\{t \geq 0 : X_t \notin D\}$ is the first exit time from $D$ by $X$. When $\mathcal{L}$ is the divergence form operator of the form (1.2), one has to run symmetric diffusion associated with $\frac{1}{2}\nabla(A\nabla)$. Observe that $X$ is in general not a semimartingale when $A$ is just measurable. Nevertheless one can still use the symmetric diffusion $X$ to solve the Dirichlet boundary value problem

$$\mathcal{L}u = 0 \qquad \text{in } D \quad \text{with} \quad u = f \qquad \text{on } \partial D$$

through a combination of Girsanov and Feynman–Kac transforms (see [8]). The Dirichlet form theory plays the role of Itô's calculus in the divergence form operator case.

In this paper, we study the Dirichlet boundary value problems for second-order elliptic operators of the following form:

$$
\begin{aligned}
(1.3) \quad \mathcal{L} &= \frac{1}{2}\nabla \cdot (A\nabla) + \mathbf{b} \cdot \nabla - \operatorname{div}(\widehat{\mathbf{b}}\cdot) + q \\
&= \frac{1}{2}\sum_{i,j=1}^{n} \frac{\partial}{\partial x_i}\left(a_{ij}(x)\frac{\partial}{\partial x_j}\right) + \sum_{i=1}^{n} b_i(x)\frac{\partial}{\partial x_i} - \text{``}\operatorname{div}(\widehat{\mathbf{b}}\cdot)\text{''} + q(x)
\end{aligned}
$$

in a bounded domain $D \subset \mathbb{R}^n$. Here $A = (a_{ij}) : \mathbb{R}^n \to \mathbb{R}^n \times \mathbb{R}^n$ is a Borel measurable, symmetric matrix-valued function that is uniformly elliptic and bounded, that is, there is a constant $\lambda \geq 1$ such that

$$(1.4) \qquad \lambda^{-1} I_{n \times n} \leq A(x) \leq \lambda I_{n \times n} \qquad \text{for every } x \in \mathbb{R}^n;$$

$\mathbf{b} = (b_1, \ldots, b_n)$ and $\widehat{\mathbf{b}} = (\widehat{b}_1, \ldots, \widehat{b}_n)$ are Borel measurable $\mathbb{R}^n$-valued functions on $\mathbb{R}^n$ and $q$ is a Borel measurable function on $\mathbb{R}^n$ such that

$$(1.5) \qquad 1_D(|\mathbf{b}|^2 + |\widehat{\mathbf{b}}|^2 + |q|) \in \mathbf{K}_n.$$

Here $\mathbf{K}_n$ denotes the space of Kato class of measures on $\mathbb{R}^n$: when $n \geq 3$, a signed measure $\mu \in \mathbf{K}_n$ if and only if

$$\lim_{r \to 0} \sup_{x \in \mathbb{R}^n} \int_{B(x,r)} |x - y|^{2-n} |\mu|(dy) = 0,$$

where $|\mu|$ denotes the total variational measure of $\mu$. Kato class $\mathbf{K}_n$ can also be defined for $n = 1, 2$; see [9] for details. A function $q$ is said to be in $\mathbf{K}_n$ if $\mu(dx) := q(x)\,dx$ is in $\mathbf{K}_n$. Clearly $L^\infty(\mathbb{R}^n) \subset \mathbf{K}_n$, and it is easy to see by using Hölder's inequality that $L^p(\mathbb{R}^n) \subset \mathbf{K}_n$ for $p > n/2$. By taking $\mathbf{b} = \widehat{\mathbf{b}} = 0 = q$ off $D$, we may and do assume in the sequel that (1.5) holds without the restriction of $D$ by $\mathbf{1}_D$.

Note that "div$(\widehat{\mathbf{b}}\cdot)$" in (1.3) is just a formal writing because the divergence really does not exist for the merely measurable vector field $\widehat{\mathbf{b}}$. It should be interpreted in the distributional sense. It is exactly due to the nondifferentiability of $\widehat{\mathbf{b}}$, all the previous known methods in solving the elliptic boundary value problems such as those in [8] and [12] ceased to work. The lower-order term div$(\widehat{\mathbf{b}}\cdot)$ cannot be handled by Girsanov transform or Feynman–Kac transform. We will show in this paper that this term in fact can be tackled by the time-reversal of Girsanov transform from the first exit time $\tau_D$ from $D$ by the symmetric diffusion $X$ associated with $\frac{1}{2}\nabla(A\nabla)$. This is the novelty of this paper. Note that time reversal of a Girsanov transform from a deterministic time was first studied in [20] for diffusions, and very recently in [4] in the context of general $m$-symmetric Markov processes. We point out that time reversal from a deterministic time in [4, 11, 20] are defined under the stationary measure $\mathbf{P}_m$. Doing time reversal from a random time $\tau_D$ involves many delicate technical issues for an effective analysis. We are able to circumvent these difficulties through a certain $h$-transform. See (1.12) below for details.

Let $(\mathcal{Q}, \mathcal{D}(\mathcal{Q}))$ be the bilinear form associated with the operator $\mathcal{L}$, where

$$\mathcal{D}(\mathcal{Q}) = \{u \in L^2(\mathbb{R}^n) : \nabla u \in L^2(\mathbb{R}^n)\} = W^{1,2}(\mathbb{R}^n)$$

and for $u, v \in W^{1,2}(\mathbb{R}^n)$

$$
\begin{aligned}
\mathcal{Q}(u,v) = {}& \frac{1}{2}\sum_{i,j=1}^n \int_{\mathbb{R}^n} a_{ij}(x)\frac{\partial u}{\partial x_i}\frac{\partial v}{\partial x_j}\,dx - \sum_{i=1}^n \int_{\mathbb{R}^n} b_i(x)\frac{\partial u}{\partial x_i}v(x)\,dx \\
& - \sum_{i=1}^n \int_{\mathbb{R}^n} \widehat{b}_i(x)\frac{\partial v}{\partial x_i}u(x)\,dx - \int_{\mathbb{R}^n} q(x)u(x)v(x)\,dx.
\end{aligned}
$$
(1.6)

For an open subset $D \subset \mathbb{R}^n$, denote by $C_c^\infty(D)$ the space of smooth functions on $D$ with compact support. The $L^2$-domain $\mathcal{D}(\mathcal{L})$ of $\mathcal{L}$ is defined to be

$$
\begin{aligned}
\{u \in W^{1,2}(\mathbb{R}^n) : \text{ there is } g \in L^2(\mathbb{R}^n) \text{ so that } \mathcal{Q}(u,v) = (-g,v)_{L^2(\mathbb{R}^n)} \\
\text{for every } v \in C_c^\infty(\mathbb{R}^n)\}
\end{aligned}
$$

and we denote $g$ by $\mathcal{L}u$. Clearly, it follows from this definition that

$$\mathcal{Q}(u,v) = (-\mathcal{L}u,v)_{L^2(\mathbb{R}^n)} \qquad \text{for } u \in \mathcal{D}(\mathcal{L}) \text{ and } v \in \mathcal{D}(\mathcal{Q}).$$

It is well known that the differential operator $\mathcal{L}$ enjoys the maximum principle if $-\operatorname{div}(\widehat{\mathbf{b}}) + q \leq 0$ in $\mathbb{R}^n$ in the following distributional sense:

$$\tag{1.7} \sum_{i=1}^{n} \int_{\mathbb{R}^n} \widehat{b}_i(x) \frac{\partial \phi}{\partial x_i} \, dx + \int_{\mathbb{R}^n} q(x)\phi(x) \, dx \leq 0$$

for all nonnegative function $\phi$ in $C_c^\infty(\mathbb{R}^n)$.

Trüdinger [26], Theorem 3.2 and Corollary 5.5, has proved the following:

THEOREM 1.1. *Assume the Markovian condition* (1.7) *holds. For every* $f \in W^{1,2}(D)$, *there exists a unique weak solution* $u \in W^{1,2}(D)$ *such that*

$$\tag{1.8} \mathcal{L}u = 0 \qquad in \ D \ with \ u - f \in W_0^{1,2}(D).$$

*Here* $W_0^{1,2}(D)$ *is the completion of* $C_c^\infty(D)$ *under the Sobolev norm*

$$\|u\|_{1,2} := \left( \int_D (u(x)^2 + |\nabla u(x)|^2) \, dx \right)^{1/2}$$

*in* $W^{1,2}(D)$. *Moreover* $u$ *is locally Hölder continuous in* $D$.

Recall that $\mathcal{L}u = 0$ in $D$ is understood in the following distributional sense:

$$\mathcal{Q}(u, \phi) = 0 \qquad \text{for every } \phi \in C_c^\infty(D).$$

We stress that condition (1.7) plays a key role in Trüdinger's approach because of the critical use of the maximum principle there.

The aim of the present paper is twofold. The first is to give a probabilistic representation for the weak solution of the Dirichlet boundary value problem (1.8). This is highly nontrivial because there is no longer a Markov process associated with the operator $\mathcal{L}$ due to the appearance of the lower-order term $\operatorname{div}(\widehat{\mathbf{b}}\cdot)$, nor can that lower-order term be handled via Girsanov transform or Feynman–Kac transform. Our idea is to use the symmetric diffusion process $X$ associated with the divergence form operator $\frac{1}{2}\nabla(A\nabla)$, the symmetric part of $\mathcal{L}$, and treat $\mathcal{L}$ as its lower-order perturbation via a combination of Girsanov and Feynman–Kac transforms and a time-reversal of Girsanov transform at the first exit time $\tau_D$ from $D$ by $X$. Based on the new probabilistic representation, our second aim is to establish the existence and uniqueness of the weak solution to problem (1.8) without the Markovian assumption (1.7). To this end, we introduce a kind of $h$-transformation which transforms the solution of the problem (1.8) to the solution of a Dirichlet boundary value problem for operators which do not involve the adjoint vector field like $\widehat{\mathbf{b}}$. The time reversal and the theory of Dirichlet forms play an essential role throughout this paper.

The remaining of the paper is organized as follows. Let $X$ be the symmetric diffusion with infinitesimal generator $\frac{1}{2}\nabla(A\nabla)$. It is well known (cf. [25]) that $X$ is a conservative Feller process on $\mathbb{R}^n$ that has Hölder continuous transition density function which admits a two-sided Aronson's Gaussian type estimate. In general, $X$ is a not a semimartingale but it admits the following Fukushima's decomposition (cf. [14]):

$$(1.9) \qquad X_t = X_0 + M_t + N_t, \qquad t \geq 0,$$

where $M = (M, \ldots, M^n)$ is a martingale additive functional (MAF) of $X$ with quadratic co-variation $\langle M^i, M^j \rangle_t = \int_0^t a_{ij}(X_s)\,ds$ and $N = (N^1, \ldots, N^n)$ is a continuous additive functional (CAF) of $X$ locally of zero quadratic variations. Note that, since $X$ has a continuous density function and each $a_{ij}$ is bounded on $\mathbb{R}^d$, by [13], Theorem 2, $M$ and $N$ can be refined to be CAFs of $X$ in the strict sense without exceptional set and (1.9) holds under $\mathbf{P}_x$ for every $x \in \mathbb{R}^n$. Without loss of generality, we work on the canonical continuous path space $C([0, \infty), \mathbb{R}^n)$ of $X$ and, for $t > 0$, denote by $r_t$ the reverse operator of $X$ from time $t$. In Section 2, we prove that, under the Markovian condition (1.7), the solution $u$ of the problem (1.8) admits the following probabilistic representation:

$$(1.10) \qquad u(x) = \mathbf{E}_x[Z_{\tau_D} f(X_{\tau_D})] \qquad \text{for } x \in D,$$

where $\tau_D := \inf\{t \geq 0 : X_t \notin D\}$ is the first exit time from $D$ by the symmetric diffusion $X$ and

$$\begin{aligned}
(1.11) \qquad Z_t = \exp\Bigg( &\int_0^t (A^{-1}\mathbf{b})(X_s)\,dM_s + \left( \int_0^t (A^{-1}\widehat{\mathbf{b}})(X_s)\,dM_s \right) \circ r_t \\
&- \frac{1}{2}\int_0^t (\mathbf{b} - \widehat{\mathbf{b}})A^{-1}(\mathbf{b} - \widehat{\mathbf{b}})^*(X_s)\,ds + \int_0^{\tau_D} q(X_s)\,ds \Bigg).
\end{aligned}$$

All the vectors in this paper are row vectors and we use $\mathbf{b}^*$ to denote the transpose of a vector $\mathbf{b}$. For two vectors $\alpha$ and $\beta$ in $\mathbb{R}^n$, we use $\alpha \cdot \beta$ or $\langle \alpha, \beta \rangle$ to denote their inner product.

Note that $Z_{\tau_D}$ is well defined under $\mathbf{P}_x$ for quasi-every $x \in \mathbb{R}^n$. This is because $t \mapsto \int_0^t (A^{-1}\widehat{\mathbf{b}})(X_s)\,dM_s$ is a square integrable martingale of $X$ in the strict sense having finite energy and so it follows from [3] and [23] that there is a continuous additive functional $L$ of $X$ having zero energy (which in general may admit an exceptional set) such that for quasi-every (q.e.) $x \in \mathbb{R}^d$, $\mathbf{P}_x$-a.s.

$$\left( \int_0^t (A^{-1}\widehat{\mathbf{b}})(X_s)\,dM_s \right) \circ r_t = -\int_0^t (A^{-1}\widehat{\mathbf{b}})(X_s)\,dM_s + L_t, \qquad t \geq 0.$$

To prove the probabilistic representation (1.10), we introduce a sequence of approximating operators $\mathcal{L}_k$ with $\widehat{\mathbf{b}}$ in the definition of $\mathcal{L}$ replaced by smooth $\widehat{\mathbf{b}}_k$. We first show that the solution $u_k$ of the Dirichlet boundary value problem for the operator $\mathcal{L}_k$ has the representation (1.10) with $\widehat{\mathbf{b}}$ replaced by smooth $\widehat{\mathbf{b}}_k$. To complete the proof, we then show that both the approximating solutions and the corresponding probabilistic expressions converge in an appropriate sense. The purpose of Section 3 is to show that the weak solution $u$ to the Dirichlet boundary value problem (1.8) is continuous up to the boundary of the domain $D$. The crucial step is to show that for every $L^p$-integrable $R^n$-valued function $f$ in a ball $B_R$ with radius $R > 0$ and $p > n$, there exists a function $v \in W_0^{1,p}(B_R)$ such that the following identity holds:

$$\left( \int_0^t f(X_s) \, dM_s \right) \circ r_t = - \int_0^t f(X_s) \, dM_s + N_t^v$$
$$\text{(1.12)}$$
$$\text{for } t < \inf\{ s : X_s \notin B_R \},$$

where $N^v$ is the zero-energy part of the Fukushima's decomposition for $v(X_t) - v(X_0)$ (see Section 3). By Sobolev embedding theorem, if we extend $v$ to $R^n$ by taking $v = 0$ on $B_R^c$, then $v$ is continuous on $R^n$. It then follows from [13], Theorem 1, that $N^v$ can be refined to be a CAF of $X$ in the strict sense. Consequently, $t \mapsto (\int_0^t f(X_s) \, dM_s) \circ r_t$ can be refined to be a CAF of $X$ in the strict sense by using the expression on the right-hand side of (1.12). In particular, $Z_{\tau_D}$ is then well defined under $\mathbf{P}_x$ for every $x \in D$. Equation (1.12) allows us to get rid of the reverse operator $r_t$ in the expression of the solution $u$ making the analysis possible. In Section 4, we establish a general existence and uniqueness result for the weak solution of the Dirichlet boundary value problem (1.8) for the operator $\mathcal{L}$ *without* the Markovian assumption (1.7). Let $Z$ be defined by (1.11). We establish as Theorem 4.4 the following important gauge theorem under suitable condition of $D$, $A$, $\mathbf{b}$, $\widehat{\mathbf{b}}$ and $q$, which is of independent interest if $\mathbf{E}_{x_0}[Z_{\tau_D}] < \infty$ for some $x_0 \in D$, then $x \mapsto \mathbf{E}_x[Z_{\tau_D}]$ is bounded between two positive constants.

We then show that if $\mathbf{E}_{x_0}[Z_{\tau_D}] < \infty$ for some $x_0 \in D$, then for every $f \in C(\partial D)$, the equation $\mathcal{L}u = 0$ in $D$ has a unique weak solution that is continuous on $\overline{D}$ such that $u = f$ on $\partial D$. Moreover, this solution can be expresses as

$$u(x) = \mathbf{E}_x[Z_{\tau_D} f(X_{\tau_D})] \qquad \text{for every } x \in D.$$

Our strategy is, by using (1.12), to reduce the Dirichlet boundary value problem $\mathcal{L}u = 0$ in $D$ with $u = f$ on $\partial D$ to a corresponding Dirichlet boundary value problem for an operator that does not have the lower-order term div($\widehat{\mathbf{b}} \cdot$). In Section 5, we consider the special case of $\mathcal{L}$ in (1.3) where $\mathbf{b} = \widehat{\mathbf{b}} = -A\nabla\rho$ where $\rho \in W^{1,2}(R^n)$ with $\nabla\rho \in L^p(R^n; dx)$ for some $p > n$.

By Sobolev embedding theorem, $\rho$ is continuous on $\mathbb{R}^n$. In this case, the quadratic form $(\mathcal{Q}, W^{1,2}(\mathbb{R}^n))$ in (1.6) takes the following form:

$$\mathcal{Q}(u,v) = \frac{1}{2}\sum_{i,j=1}^{n}\int_{\mathbb{R}^n} a_{ij}(x)\frac{\partial u}{\partial x_i}\frac{\partial v}{\partial x_j}\,dx + \frac{1}{2}\sum_{i,j=1}^{n}\int_{\mathbb{R}^n} a_{ij}(x)\frac{\partial(uv)}{\partial x_i}\frac{\partial \rho}{\partial x_j}\,dx$$

$$+ \int_{\mathbb{R}^n} u(x)v(x)q(x)\,dx$$

for $u,v \in W^{1,2}(\mathbb{R}^n)$. For this case, we can establish a stronger result without additional condition on the diffusion matrix $A$. Let $X$ be the symmetric diffusion with infinitesimal generator $\frac{1}{2}\nabla(A\nabla)$, and recall the Fukushima's decomposition in the strict sense (cf. [13], Theorem 1):

$$\rho(X_t) - \rho(X_0) = M_t^\rho + N_t^\rho, \qquad t \geq 0,$$

where $M^\rho$ is the MAF of $X$ in the strict sense having finite energy and $N^\rho$ is the continuous additive functional of $X$ in the strict sense having zero energy. Define

$$Z_t = \exp\left(N_t^\rho + \int_0^t q(X_s)\,ds\right), \qquad t \geq 0.$$

We prove in Theorem 5.1 that if $D$ is a bounded Lipschitz domain and if $\mathbf{E}_{x_0}[Z_{\tau_D}] < \infty$ for some $x_0 \in D$, then $x \mapsto \mathbf{E}_x[Z_{\tau_D}]$ is bounded between two positive constants. Moreover, assuming that $\mathbf{E}_{x_0}[Z_{\tau_D}] < \infty$ for some $x_0 \in D$, we show that for every $f \in C(\partial D)$,

$$u(x) := \mathbf{E}_x[Z_{\tau_D}f(X_{\tau_D})], \qquad x \in D,$$

gives the unique weak solution of $\mathcal{L}u = 0$ in $D$ that is continuous on $\overline{D}$ with $u = f$ on $\partial D$.

In this paper, we use ":=" as a way of definition. A statement is said to hold quasi-everywhere (q.e.) on some set $A \subset \mathbb{R}^n$ if there is an exceptional set $N$ of zero capacity so that the statement holds on $A \setminus N$. For the general theory of Dirichlet forms and Markov processes and their terminology, we refer readers to [14] and [21] .

**2. Probabilistic representation.** In this section, we will give a probabilistic representation of the weak solutions of Dirichlet boundary value problems. Consider the following regular Dirichlet form on $\mathbb{R}^n$:

$$\mathcal{E}(u,v) = \frac{1}{2}\sum_{i,j=1}^{n}\int_{\mathbb{R}^n} a_{ij}(x)\frac{\partial u}{\partial x_i}\frac{\partial v}{\partial x_j}\,dx,$$

(2.1)

$$\mathcal{D}(\mathcal{E}) = W^{1,2}(\mathbb{R}^n).$$

It is well known that there is a symmetric conservative diffusion process $X = \{X_t, \theta_t, r_t, P_x, x \in \mathbb{R}^n\}$ associated with it. Since $X$ has Hölder continuous transition density function with respect to the Lebesgue measure on $\mathbb{R}^n$ (cf. [25]), $X$ can be modified to start from every point in $\mathbb{R}^n$. Without loss of generality, we may and do assume that $X$ is defined on the canonical sample space $\Omega = C([0, \infty) \to \mathbb{R}^n)$ on which the time-shift operators $\{\theta_t, t \geq 0\}$ and time-reversal operators $\{r_t, t > 0\}$ are well defined: for $t > 0$,

$$(2.2) \qquad \theta_t(\omega)(s) = \omega(t + s) \qquad \text{for } s \geq 0$$

and

$$(2.3) \qquad r_t(\omega)(s) := \begin{cases} \omega(t - s), & \text{if } 0 \leq s \leq t, \\ \omega(0), & \text{if } s \geq t. \end{cases}$$

Let $\{\mathcal{F}_t, t \geq 0\}$ be the minimal augmented filtration generated by the diffusion process $X$. For every $u \in D(\mathcal{E})$, the following Fukushima's decomposition holds: for q.e. $x \in \mathbb{R}^n$,

$$(2.4) \qquad u(X_t) - u(X_0) = M_t^u + N_t^u, \qquad \mathbf{P}_x\text{-a.s.},$$

where $M^u$ a continuous MAF of $X$ having finite energy and $N_t^u$ is a CAF of $X$ having zero energy. Note that the MAF $M^u$ and the CAF $N^u$ typically admit an exceptional set $N$ of zero capacity in their definition. However since $X$ has a Hölder continuous transition density function, by [13], Theorem 1, both $M^u$, $N^u$ and the above Fukushima decomposition (2.4) can be strengthened to admit no exceptional set [so in particular, (2.4) holds for every $x \in \mathbb{R}^n$] if $u$ is continuous and the energy measure for $M^u$

$$\mu_{\langle u \rangle}(dx) := \sum_{i,j=1}^n a_{i,j}(x) \frac{\partial u(x)}{\partial x_i} \frac{\partial u(x)}{\partial x_j} \, dx$$

is a smooth measure in the strict sense. The latter means that there is an increasing sequence of finely open sets $\{D_k, k \geq 1\}$ so that $\bigcup_{k=1}^\infty D_k = \mathbb{R}^n$, $1_{D_k} \mu_{\langle u \rangle}$ is a finite Borel measure and $G_1(1_{D_k} \mu_{\langle u \rangle})$ is bounded for every $k \geq 1$. In the sequel we call an additive functional strict if it admits no exceptional set and call the refined decomposition of (2.4) without exceptional set a strict Fukushima decomposition. In fact by [13], Theorem 2, both $M^u$ and $N^u$ can be taken to be strict AFs of $X$ and the strict Fukushima decomposition holds for every continuous function $u$ that is locally in $\mathcal{F}$ such that $\mu_{\langle u \rangle}$ is a smooth measure in the strict sense. Applying the above to coordinate functions $f_j(x) := x_j$ for $j = 1, \dots, n$, we have

$$(2.5) \qquad X_t = x + M_t + N_t, \qquad \mathbf{P}_x\text{-a.s.}$$

for every $x \in \mathbb{R}^n$, where $M_t = (M_t^1, \dots, M_t^n)$ is a continuous local MAF of $X$ in the strict sense with

$$\langle M^i, M^j \rangle_t = \int_0^t a_{ij}(X_s) \, ds$$

and $N_t$ is a CAF of $X$ locally of zero energy in the strict sense. In particular, there is a Brownian motion $B = (B^1, \ldots, B^n)$, which is a martingale AF of $X$ in the strict sense, such that

$$M = \int_0^t \sigma(X_s) \, dB_s, \qquad t \geq 0,$$

where $\sigma(x)$ is the positive definite symmetric square root of the matrix $A(x)$. Note that $M$ is a MAF of $X$ in the strict sense.

Here is the first representation result.

LEMMA 2.1. *Assume condition (1.7) holds and that $\widehat{\mathbf{b}}$ is $C^1$-smooth. Then for every $f \in W^{1,2}(D) \cap C(\overline{D})$, the unique weak solution $u$ to (1.8) admits the following representation: for $x \in D$*

$$
\begin{aligned}
(2.6) \quad u(x) = \mathbf{E}_x \Big[ f(X_{\tau_D}) \exp \Big\{ & \int_0^{\tau_D} (A^{-1} \mathbf{b})(X_s) \, dM_s \\
& + \Big( \int_0^{\tau_D} (A^{-1} \widehat{\mathbf{b}})(X_s) \, dM_s \Big) \circ r_{\tau_D} \\
& - \frac{1}{2} \int_0^{\tau_D} (\mathbf{b} - \widehat{\mathbf{b}}) A^{-1} (\mathbf{b} - \widehat{\mathbf{b}})^* (X_s) \, ds \\
& + \int_0^{\tau_D} q(X_s) \, ds \Big\} \Big].
\end{aligned}
$$

*Moreover, $u \in C(\overline{D})$.*

PROOF. First we note that by [20], (46), we for a.e. $x \in D$, $\mathbf{P}_x$-a.s. have

$$
\begin{aligned}
(2.7) \quad & \Big( \int_0^{\tau_D} (A^{-1} \widehat{\mathbf{b}})(X_s) \, dM_s \Big) \circ r_{\tau_D} \\
& = - \int_0^{\tau_D} (A^{-1} \widehat{\mathbf{b}})(X_s) \, dM_s - \int_0^{\tau_D} \operatorname{div}(\widehat{\mathbf{b}})(X_s) \, ds.
\end{aligned}
$$

Note that since $\widehat{\mathbf{b}} \in C^1(\mathbb{R}^n)$, both $t \mapsto \int_0^t (A^{-1} \widehat{\mathbf{b}})(X_s) \, dM_s$ and $t \mapsto \int_0^t \operatorname{div}(\widehat{\mathbf{b}}) \times (X_s) \, ds$ are continuous AFs of $X$ in the strict sense. Therefore $(\int_0^{\tau_D} (A^{-1} \widehat{\mathbf{b}}) \times (X_s) \, dM_s) \circ r_{\tau_D}$ can be refined using the right-hand side of (2.7) so that it is well defined under $\mathbf{P}_x$ for every $x \in D$. In particular, the right-hand side of (2.6) is well defined for every $x \in D$.

Under our assumptions, the operator $\mathcal{L}$ can be written as

$$
\begin{aligned}
\mathcal{L} = \frac{1}{2} \sum_{i,j=1}^n & \frac{\partial}{\partial x_i} \Big( a_{ij}(x) \frac{\partial}{\partial x_j} \Big) \\
& + \sum_{i=1}^n (b_i(x) - \widehat{b}_i(x)) \frac{\partial}{\partial x_i} - \operatorname{div} \widehat{\mathbf{b}}(x) + q(x).
\end{aligned}
$$

Define a family of measures $\{\mathbf{Q}_x, x \in \mathbb{R}^n\}$ on $\mathcal{F}_\infty$ by

$$(2.8) \qquad \frac{d\mathbf{Q}_x}{d\mathbf{P}_x}\bigg|_{\mathcal{F}_t} = H_t,$$

where

$$H_t = \exp\left(\int_0^t A^{-1}(\mathbf{b}-\widehat{\mathbf{b}})(X_s)\,dM_s - \frac{1}{2}\int_0^t (\mathbf{b}-\widehat{\mathbf{b}})A^{-1}(\mathbf{b}-\widehat{\mathbf{b}})^*(X_s)\,ds\right).$$
(2.9)

It is known (cf. [8]) that under measure $\{\mathbf{Q}_x, x \in \mathbb{R}^n\}$, $X$ is a diffusion process on $\mathbb{R}^n$ having infinitesimal generator

$$\mathcal{L}_0 = \frac{1}{2}\sum_{i,j=1}^n \frac{\partial}{\partial x_i}\left(a_{ij}(x)\frac{\partial}{\partial x_j}\right) + \sum_{i=1}^n (b_i(x) - \widehat{b}_i(x))\frac{\partial}{\partial x_i}.$$

By Theorem 5.11 in [8], when $f \in W^{1,2}(D) \cap C(\overline{D})$, the unique weak solution $u$ to (1.8) is continuous on $\overline{D}$ and has the following probabilistic representation

$$u(x) = \mathbf{E}_x^Q\left[f(X_{\tau_D})\exp\left(\int_0^{\tau_D}(-\operatorname{div}(\widehat{\mathbf{b}})(X_s) + q(X_s))\,ds\right)\right], \qquad x \in D,$$

where $\mathbf{E}_x^Q$ stands for the expectation with respect to the measure $\mathbf{Q}_x$. Since $\{H_{t\wedge\tau_D}, t \geq 0\}$ is a uniformly integrable martingale under $\mathbf{P}_x$ for every $x \in D$, due to the fact that $|\mathbf{b}-\widehat{\mathbf{b}}|^2 \in \mathbf{K}_n$ (see [6], page 746), we have

$$u(x) = \mathbf{E}_x\bigg[f(X_{\tau_D})\exp\left(\int_0^{\tau_D}(A^{-1}(\mathbf{b}-\widehat{\mathbf{b}}))(X_s)\,dM_s\right.$$
$$-\frac{1}{2}\int_0^{\tau_D}(\mathbf{b}-\widehat{\mathbf{b}})A^{-1}(\mathbf{b}-\widehat{\mathbf{b}})^*(X_s)\,ds\bigg)$$
$$\times\exp\left(\int_0^{\tau_D}(-\operatorname{div}(\widehat{\mathbf{b}})(X_s) + q(X_s))\,ds\right)\bigg].$$

This together with (2.7) implies (2.6). $\quad\square$

Put

$$J(x) = \begin{cases} c_0\exp\left(-\dfrac{1}{1-|x|^2}\right), & \text{if } |x| < 1, \\ 0, & \text{if } |x| \geq 1, \end{cases}$$

where $c_0 > 0$ is a normalizing constant so that

$$\int_{\mathbb{R}^n} J(x)\,dx = 1.$$

For any positive integer $k \geq 1$, we set

$$J_k(x) := k^n J(kx) \qquad \text{for } x \in \mathbb{R}^n$$

and

$$\widehat{\mathbf{b}}_k(x) = J_k * \widehat{\mathbf{b}}(x) := \int_{\mathbb{R}^n} \widehat{\mathbf{b}}(y) J_k(x - y) \, dy,$$

(2.10)

$$q_k(x) = J_k * q(x) := \int_{\mathbb{R}^n} q(y) J_k(x - y) \, dy.$$

Recall that we assume that $\mathbf{b}$, $\widehat{\mathbf{b}}$ and $q$ are set to be zero off $D$, and $|\mathbf{b}|^2 + |\widehat{\mathbf{b}}|^2 + |q| \in \mathbf{K}_n$.

LEMMA 2.2. (i) $\widehat{\mathbf{b}}_k \to \widehat{\mathbf{b}}$ in $L^2_{\mathrm{loc}}(\mathbb{R}^n)$ and $q_k \to q$ in $L^1_{\mathrm{loc}}(\mathbb{R}^n)$ as $k \to \infty$.

(ii) For every nonnegative function $\phi$ in $C_c^\infty(\mathbb{R}^n)$ and $k \geq 1$,

$$\sum_{i=1}^n \int_{\mathbb{R}^n} \widehat{b}_{ki}(x) \frac{\partial \phi}{\partial x_i} \, dx + \int_{\mathbb{R}^n} q_k(x) \phi(x) \, dx \leq 0.$$

(iii) Assume $n \geq 3$. Then

$$\lim_{r \to 0} \sup_{\substack{k \geq 1 \\ x \in R^n}} \int_{\{y \in \mathbb{R}^n : |x-y| \leq r\}} \frac{|\widehat{\mathbf{b}}_k(y)|^2 + |q_k(y)|}{|x - y|^{n-2}} \, dy = 0.$$

(iv) Assume $n \geq 3$. Then

$$\lim_{k \to \infty} \sup_{x \in \mathbb{R}^n} \int_D \frac{|\widehat{\mathbf{b}}_k(y) - \widehat{\mathbf{b}}(y)|^2 + |q_k(y) - q(y)|}{|x - y|^{n-2}} \, dy = 0.$$

PROOF. (i) follows easily from the fact that $\widehat{\mathbf{b}} \in L^2_{\mathrm{loc}}(\mathbb{R}^n)$ and $q \in L^1_{\mathrm{loc}}(\mathbb{R}^n)$. For a nonnegative function $\phi$ in $C_c^\infty(\mathbb{R}^n)$ and $z \in \mathbb{R}^n$, put $\phi_z(x) = \phi(x + z)$. Then it is easy to see that for $k \geq 1$,

$$\sum_{i=1}^n \int_{\mathbb{R}^n} \widehat{b}_{ki}(x) \frac{\partial \phi}{\partial x_i} \, dx + \int_{\mathbb{R}^n} q_k(x) \phi(x) \, dx$$

$$= \int_{\mathbb{R}^n} J_k(z) \left( \sum_{i=1}^n \int_{\mathbb{R}^n} \widehat{b}_i(x) \frac{\partial \phi_z}{\partial x_i} \, dx + \int_{\mathbb{R}^n} q(x) \phi_z(x) \, dx \right) dz$$

and so (ii) follows from (1.7). Since $|\widehat{\mathbf{b}}|^2 \in \mathbf{K}_n$, (iii) is a consequence of the following inequality:

$$\sup_{\substack{x \in \mathbb{R}^n, |x-y| \leq r \\ k \geq 1}} \int_{|x-y| \leq r} \frac{|\widehat{\mathbf{b}}_k|^2(y)}{|x - y|^{n-2}} \, dy \leq \int_{|x-y| \leq r} \frac{\int_{\mathbb{R}^n} J_k(z) |\widehat{\mathbf{b}}|^2(y-z) \, dz}{|x - y|^{n-2}} \, dy$$

$$= \int_{\mathbb{R}^n} J_k(z) \left( \int_{|x-y| \leq r} \frac{|\widehat{\mathbf{b}}|^2(y-z)}{|x - y|^{n-2}} \, dy \right) dz$$

$$= \int_{\mathbb{R}^n} J_k(z) \left( \int_{|x-z-y| \le r} \frac{|\widehat{\mathbf{b}}|^2(y)}{|x-z-y|^{n-2}} \, dy \right) dz$$

$$\le \sup_{x \in \mathbb{R}^n} \int_{|x-y| \le r} \frac{|\widehat{\mathbf{b}}|^2(y)}{|x-y|^{n-2}} \, dy,$$

where we used the fact $\int_{\mathbb{R}^n} J_k(z) \, dz = 1$. The proof for $q$ is similar. To prove (iv) we observe from the proof of (iii) that for any $r > 0$,

$$\sup_{x \in \mathbb{R}^n} \int_{\mathbb{R}^n} \frac{|\widehat{\mathbf{b}}_k - \widehat{\mathbf{b}}|^2(y)}{|x-y|^{n-2}} \, dy$$

$$\le 2 \sup_{x \in \mathbb{R}^n} \int_{|x-y| \le r} \frac{|\widehat{\mathbf{b}}|^2(y)}{|x-y|^{n-2}} \, dy + \frac{1}{r^{n-2}} \int_{\mathbb{R}^n} |\widehat{\mathbf{b}}_k - \widehat{\mathbf{b}}|^2(y) \, dy.$$

Using (i), this implies that

$$\lim_{k \to \infty} \sup_{x \in \mathbb{R}^n} \int_{\mathbb{R}^n} \frac{|\widehat{\mathbf{b}}_k - \widehat{\mathbf{b}}|^2(y)}{|x-y|^{n-2}} \, dy \le 2 \sup_{x \in \mathbb{R}^n} \int_{|x-y| \le r} \frac{|\widehat{\mathbf{b}}|^2(y)}{|x-y|^{n-2}} \, dy$$

for any $r > 0$. Letting $r \to 0$ we get (iv) for $\widehat{\mathbf{b}}_k - \widehat{\mathbf{b}}$. The proof for $|q_k - q|$ goes in a similar way. $\square$

For integer $k \ge 1$, define

$$
\begin{aligned}
(2.11) \quad \mathcal{L}_k := & \frac{1}{2} \sum_{i,j=1}^{n} \frac{\partial}{\partial x_i} \left( a_{ij}(x) \frac{\partial}{\partial x_j} \right) \\
& + (\mathbf{b}(x) - \widehat{\mathbf{b}}_k(x)) \cdot \nabla - \operatorname{div} \widehat{\mathbf{b}}_k(x) + q_k(x),
\end{aligned}
$$

where $\widehat{\mathbf{b}}_k, q_k$ are the functions defined by (2.10). Denote by $\mathcal{Q}_k(\cdot, \cdot)$ the quadratic form associated with $\mathcal{L}_k$.

Now we can drop the assumption of $\widehat{\mathbf{b}}$ being $C^1$ from Lemma 2.1 by using the smooth approximation $\widehat{\mathbf{b}}_k$ for $\widehat{\mathbf{b}}$.

THEOREM 2.3. *Suppose that condition (1.7) holds and $f \in W^{1,2}(D) \cap C(\overline{D})$. Then the unique weak solution $u$ of (1.8) has the following probabilistic representation: for q.e. $x \in D$,*

$$
\begin{aligned}
(2.12) \quad u(x) = \mathbf{E}_x \bigg[ & f(X_{\tau_D}) \exp \bigg( \int_0^{\tau_D} (A^{-1} \mathbf{b})(X_s) \, dM_s \\
& + \bigg( \int_0^{\tau_D} (A^{-1} \widehat{\mathbf{b}})(X_s) \, dM_s \bigg) \circ r_{\tau_D} \\
& - \frac{1}{2} \int_0^{\tau_D} (\mathbf{b} - \widehat{\mathbf{b}}) A^{-1} (\mathbf{b} - \widehat{\mathbf{b}})^*(X_s) \, ds \\
& + \int_0^{\tau_D} q(X_s) \, ds \bigg) \bigg].
\end{aligned}
$$

PROOF.   For simplicity, we assume $n \geq 3$. [The case of $n = 1$ and $n = 2$ can be handled similarly with a corresponding version for Lemma 2.2(iii) and (iv).] Recall the differential operator $\mathcal{L}_k$ defined by (2.11). Let $u_k$ denote the unique weak solution of the following Dirichlet boundary value problem:

$$\mathcal{L}_k u_k = 0 \text{ in } D \qquad \text{with } u_k - f \in W_0^{1,2}(D).$$

Since $\widehat{\mathbf{b}}_k$ is smooth, it follows from Lemmas 2.2(ii) and 2.1 that

$$
\begin{aligned}
u_k(x) = \mathbf{E}_x \Bigg[ f(X_{\tau_D}) \exp\Bigg( & \int_0^{\tau_D} (A^{-1}\mathbf{b})(X_s)\, dM_s \\
& + \left( \int_0^{\tau_D} (A^{-1}\widehat{\mathbf{b}}_k)(X_s)\, dM_s \right) \circ r_{\tau_D} \\
& - \frac{1}{2} \int_0^{\tau_D} (\mathbf{b} - \widehat{\mathbf{b}}_k) A^{-1} (\mathbf{b} - \widehat{\mathbf{b}}_k)^*(X_s)\, ds \\
& \qquad\qquad + \int_0^{\tau_D} q_k(X_s)\, ds \Bigg) \Bigg].
\end{aligned}
$$

Let $v$ denote the right-hand side of (2.12). We will show that $\lim_{k \to \infty} u_k(x) = v(x)$. To this end, put

$$
\begin{aligned}
Z_k := \exp\Bigg( & \int_0^{\tau_D} (A^{-1}\mathbf{b})(X_s)\, dM_s + \left( \int_0^{\tau_D} (A^{-1}\widehat{\mathbf{b}}_k)(X_s)\, dM_s \right) \circ r_{\tau_D} \\
& - \frac{1}{2} \int_0^{\tau_D} (\mathbf{b} - \widehat{\mathbf{b}}_k) A^{-1} (\mathbf{b} - \widehat{\mathbf{b}}_k)^*(X_s)\, ds + \int_0^{\tau_D} q_k(X_s)\, ds \Bigg)
\end{aligned}
$$

and

$$
\begin{aligned}
Z := \exp\Bigg( & \int_0^{\tau_D} (A^{-1}\mathbf{b})(X_s)\, dM_s + \left( \int_0^{\tau_D} (A^{-1}\widehat{\mathbf{b}})(X_s)\, dM_s \right) \circ r_{\tau_D} \\
& - \frac{1}{2} \int_0^{\tau_D} (\mathbf{b} - \widehat{\mathbf{b}}) A^{-1} (\mathbf{b} - \widehat{\mathbf{b}})^*(X_s)\, ds + \int_0^{\tau_D} q(X_s)\, ds \Bigg).
\end{aligned}
$$

We first prove that $Z_k \to Z$ in probability as $k \to \infty$. It is clear that

$$\int_0^{\tau_D} (\mathbf{b} - \widehat{\mathbf{b}}_k) A^{-1} (\mathbf{b} - \widehat{\mathbf{b}}_k)^*(X_s)\, ds + \int_0^{\tau_D} q_k(X_s)\, ds$$

converges in probability under $\mathbf{P}_x$ for every $x \in D$ to

$$\int_0^{\tau_D} (\mathbf{b} - \widehat{\mathbf{b}}) A^{-1} (\mathbf{b} - \widehat{\mathbf{b}})^*(X_s)\, ds + \int_0^{\tau_D} q(X_s)\, ds.$$

Thus, it is sufficient to show that

$$(2.13) \quad \left( \int_0^{\tau_D} (A^{-1}\widehat{\mathbf{b}}_k)(X_s)\, dM_s \right) \circ r_{\tau_D} \to \left( \int_0^{\tau_D} (A^{-1}\widehat{\mathbf{b}})(X_s)\, dM_s \right) \circ r_{\tau_D}$$

in probability under $\mathbf{P}_x$ as $k \to \infty$ for q.e. $x \in D$. Define

$$\widehat{M}_t^k := \int_0^t (A^{-1} \widehat{\mathbf{b}}_k)(X_s) \, dM_s, \qquad t \geq 0$$

and

$$\widehat{M}_t := \int_0^t (A^{-1} \widehat{\mathbf{b}})(X_s) \, dM_s, \qquad t \geq 0,$$

which are MAFs of $X$ in the strict sense of finite energy (recall that we assumed $\mathbf{b} = \widehat{\mathbf{b}} = \mathbf{0}$ off $D$). It follows from [3] and [23] that there are continuous processes $N_t^k$ and $N_t$ of zero energy such that

$$\widehat{M}_t^k \circ r_t = -\widehat{M}_t^k + N_t^k \quad \text{and} \quad \widehat{M}_t \circ r_t = -\widehat{M}_t + N_t.$$

Moreover, since $\widehat{M}^k \to \widehat{M}$ as $k \to \infty$ with respect to the energy norm in the martingale space, for every subsequence $\{n_k\}$, there is a sub-subsequence $\{n_{k_j}\}$ so that $N_t^{n_{k_j}}$ converges to $N_t$ uniformly on compact intervals $\mathbf{P}_x$-a.s. for q.e. $x \in D$ (cf. [14]). Thus, for any $T > 0$, on $\{\omega : \tau_D(\omega) \leq T\}$, it holds that

$$\left| \left( \int_0^{\tau_D} (A^{-1} \widehat{\mathbf{b}}_k)(X_s) \, dM_s \right) \circ r_{\tau_D} - \left( \int_0^{\tau_D} (A^{-1} \widehat{\mathbf{b}})(X_s) \, dM_s \right) \circ r_{\tau_D} \right|$$

$$\leq \sup_{0 \leq t \leq T} |\widehat{M}_t^k \circ r_t - \widehat{M}_t \circ r_t|$$

$$\leq \sup_{0 \leq t \leq T} |M_t^k - \widehat{M}_t| + \sup_{0 \leq u \leq T} |N_t^k - N_t|.$$

This proves (2.13). So, to show that $u_n \to v$, it suffices to prove that the family $\{Z_n, n \geq 1\}$ is uniformly integrable under $\mathbf{P}_x$ for q.e. $x \in D$. In view of Lemmas 2.1, (2.7) and 2.2(ii), we have

$$Z_k \leq L_k := \exp\left( \int_0^{\tau_D} (A^{-1}(\mathbf{b} - \widehat{\mathbf{b}}_k))(X_s) \, dM_s \right.$$

$$\left. - \frac{1}{2} \int_0^{\tau_D} (\mathbf{b} - \widehat{\mathbf{b}}_k) A^{-1} (\mathbf{b} - \widehat{\mathbf{b}}_k)^*(X_s) \, ds \right).$$

Define $L$ in the same way as $L_k$ but with $\widehat{\mathbf{b}}$ in place of $\widehat{\mathbf{b}}_k$. Then

$$\int_\Omega |L_k - L| \, d\mathbf{P}_x = \int_{\{L_k > L\}} (L_k - L) \, d\mathbf{P}_x + \int_{\{L_k \leq L\}} (L - L_k) \, d\mathbf{P}_x$$

$$= \int_\Omega L_k \, d\mathbf{P}_x - \int_\Omega L \, d\mathbf{P}_x + 2 \int_{\{L_k \leq L\}} (L - L_k) \, d\mathbf{P}_x$$

$$= 2 \int_{\{L_k \leq L\}} (L - L_k) \, d\mathbf{P}_x \to 0 \qquad \text{as } k \to \infty.$$

Here we have used the fact that $\mathbf{E}_x[L_k] = 1 = \mathbf{E}_x[L]$, which is a consequence of the Kato class assumption on $|\mathbf{b}|^2 + |\hat{\mathbf{b}}|^2$ (cf. [2]). This particularly implies that $\{L_k, k \geq 1\}$ is uniformly integrable under $\mathbf{P}_x$ for every $x \in D$, so is $\{Z_k, k \geq 1\}$.

To show that the weak solution $u$ of (1.8) is equal to $v$, by the uniqueness, it suffices to show that $v$ is a weak solution to (1.8). By Theorem 3.2 (and its proof) of Trüdinger [26], there is a constant $C > 0$, independent of $k$, such that

$$\|u_k\|_{1,2} \leq C\|f\|_{1,2} \qquad \text{for every } k \geq 1.$$

By taking a subsequence if necessary, we may assume that $u_n$ converges weakly to some $v_1$ in $W^{1,2}(D)$ and that its Cesaro mean $\{k^{-1}\sum_{j=1}^k u_j, k \geq 1\}$ converges to some $v_2$ in $(W^{1,2}(D), \|\cdot\|_{1,2})$. Clearly $v_1 = v_2 = v$. Moreover, since $u_k - f \in W_0^{1,2}(D)$, we have $v - f \in W_0^{1,2}(D)$.

As $\mathcal{Q}_k(u_k, \phi) = 0$ for any $\phi \in C_c^\infty(D)$, if we can show that for any $\phi \in C_c^\infty(D)$, $\mathcal{Q}(v, \phi) = \lim_{k\to\infty} \mathcal{Q}_k(u_k, \phi)$, then $v \in W^{1,2}(D)$ is a weak solution for $\mathcal{L}v = 0$ in $D$ with $v - f \in W_0^{1,2}(D)$.

Observe that

$$\mathcal{Q}_k(u_k, \phi) = \frac{1}{2}\sum_{i,j=1}^n \int_{\mathbb{R}^n} a_{ij}(x)\frac{\partial u_k}{\partial x_i}\frac{\partial \phi}{\partial x_j}\,dx - \sum_{i=1}^n \int_{\mathbb{R}^n} b_i(x)\frac{\partial u_k}{\partial x_i}\phi\,dx$$

$$- \sum_{i=1}^n \int_{\mathbb{R}^n} \hat{b}_{ki}(x)\frac{\partial \phi}{\partial x_i}u_k(x)\,dx - \int_{\mathbb{R}^n} q_k(x)u_k(x)\phi\,dx.$$

Obviously

$$\lim_{k\to\infty}\frac{1}{2}\sum_{i,j=1}^n \int_{\mathbb{R}^n} a_{ij}(x)\frac{\partial u_k}{\partial x_i}\frac{\partial \phi}{\partial x_j}\,dx = \frac{1}{2}\sum_{i,j=1}^n \int_{\mathbb{R}^n} a_{ij}(x)\frac{\partial v}{\partial x_i}\frac{\partial \phi}{\partial x_j}\,dx$$

and

$$\lim_{k\to\infty}\sum_{i=1}^n \int_{\mathbb{R}^n} b_i(x)\frac{\partial u_k(x)}{\partial x_i}\phi(x)\,dx = \sum_{i=1}^n \int_{\mathbb{R}^n} b_i(x)\frac{\partial v(x)}{\partial x_i}\phi(x)\,dx.$$

On the other hand, for $\phi \in C_c^\infty(D)$, by Lemma 2.2 we have

$$\lim_{k\to\infty}\sum_{i=1}^n \int_{\mathbb{R}^n} \hat{b}_{ki}(x)\frac{\partial \phi(x)}{\partial x_i}u_k(x)\,dx = \sum_{i=1}^n \int_{\mathbb{R}^n} \hat{b}_i(x)\frac{\partial \phi(x)}{\partial x_i}v(x)\,dx$$

and

$$\lim_{k\to\infty}\int_{\mathbb{R}^n} q_k(x)u_k(x)\phi(x)\,dx = \int_{\mathbb{R}^n} q(x)v(x)\phi(x)\,dx.$$

This proves that

$$\mathcal{Q}(v, \phi) = \lim_{k\to\infty}\mathcal{Q}_k(u_k, \phi) = 0 \qquad \text{for every } \phi \in C_c^\infty(D)$$

and so $\mathcal{L}v = 0$ in $D$ in the distributional sense. This proves the theorem. $\quad\square$

**3. Continuity at the boundary.** In this section, we study the regularity of the solution of the Dirichlet boundary value problems (1.8) at the boundary of the domain. First, we prepare two useful lemmas. The next result is due to Meyers [22], Theorem 1.

LEMMA 3.1. *For every $x_0 \in \mathbb{R}^n$, $R > 0$ and $p > n$, there is a constant $\varepsilon \in (0, 1)$, depending only on $n, R$ and $p$, such that if*

$$(3.1) \quad (1 - \varepsilon)I_{n \times n} \leq A(x) \leq I_{n \times n} \qquad \text{for a.e. } x \in B_R := B(x_0, R),$$

*then $\frac{1}{2}\nabla(A\nabla u) = \operatorname{div} \mathbf{f}$ in $B_R$ has a unique weak solution in $W_0^{1,p}(B_R)$ for every $\mathbf{f} = (f_1, \ldots, f_n) \in L^2(B_R; dx)$. Moreover, there is a constant $c > 0$ independent of $\mathbf{f}$ such that*

$$\|\nabla u\|_{L^p(B_R; dx)} \leq c\|\mathbf{f}\|_{L^p(B_R; dx)}.$$

Observe that since $u \in W_0^{1,p}(D) \subset \mathbb{R}^n$ with $p > n$, by the classical Sobolev embedding theorem (see, e.g., [16], Theorem 7.10) $u \in C(\overline{D})$ if we take $u = 0$ on $D^c$. Recall that $D$ is a bounded domain in $\mathbb{R}^n$. Select $x_0 \in \mathbb{R}^n$ and $R > 0$ so that $B(x_0, R) \supset D$. For simplicity, denote $B(x_0, R)$ by $B_R$. Let $X^{B_R}$ denote the symmetric diffusion in $B_R$ associated with the infinitesimal generator $\frac{1}{2}\nabla(A\nabla)$. The time-reversal operator for $X^{B_R}$ will still be denoted as $r_t$.

For a MAF $\widehat{M}$ of $X^{B_R}$ of finite energy, let

$$\Gamma(\widehat{M})_t = -\tfrac{1}{2}(\widehat{M}_t + \widehat{M}_t \circ r_t) \qquad \text{for } t < \tau_{B_R}.$$

According to the representation theorem of martingale additive functionals in [14], there is a measurable function $F : B_R \to \mathbb{R}^n$ such that

$$\widehat{M}_t = \int_0^t F(X_s) \, dM_s \qquad \text{for } t < \tau_{B_R}.$$

We have the following result.

LEMMA 3.2. *Suppose that $p > n$ and the diffusion matrix $A$ satisfies the condition (3.1). Then for every $L^p$-integrable $\mathbb{R}^n$-valued function $F : B_R \to \mathbb{R}^n$, there exists a function $v \in W_0^{1,p}(B_R) \subset W_0^{1,2}(B_R)$ with $\|\nabla v\|_{L^p(B_R; dx)} \leq c\|F\|_{L^p(B_R; dx)}$ such that*

$$(3.2) \quad \Gamma(\widehat{M})_t = \Gamma(M^v)_t = N_t^v \qquad \text{for } t < \tau_{B_R}.$$

*Moreover, if we extend $v$ to $\mathbb{R}^n$ by taking $v = 0$ on $B_R^c$, then $v \in C(\mathbb{R}^n)$.*

PROOF. By [10], Corollary 3.2, the following orthogonal decomposition holds with respect to the inner product induced by the energy norm.

$$(3.3) \quad \widehat{M}_t = M_t^v + K_t \qquad \text{for } t < \tau_{B_R},$$

where $v$ is an element in $W_0^{1,2}(B_R)$, $K$ is a MAF of $X^{B_R}$ satisfying $\Gamma(K) = 0$ and $\mu_{\langle K, M^\phi \rangle}(B_R) = 0$ for all $\phi \in D(\mathcal{E})$. Here $\mu_{\langle K, M^\phi \rangle}$ is the signed Revuz measure of CAF $\langle K, M^\phi \rangle$ of $X$ of finite variation. Hence,

$$\Gamma(\widehat{M})_t = \Gamma(M^v)_t = N_t^v \qquad \text{for } t < \tau_{B_R}.$$

By the representation of MAFs, we have

$$M_t^v = \int_0^t \nabla v(X_s) \, dM_s \quad \text{and} \quad K_t = \int_0^t G(X_s) \, dM_s \qquad \text{for } t < \tau_{B_R},$$

for some measurable vector field $G = (G_1(x), \ldots, G_n(x)) : B_R \to \mathbb{R}^n$. Since $\mu_{\langle K, M^\phi \rangle}(B_R) = 0$ for all $\phi \in W_0^{1,2}(B_R)$, we have

$$\int_D \sum_{i,j=1}^n a_{ij}(x) \frac{\partial \phi}{\partial x_i} G_j(x) \, dx = 0$$

for all $\phi \in C_c^1(B_R)$. This says that $\operatorname{div}(AG) = 0$ in $B_R$. Note that

$$F(x) = \nabla v(x) + G(x).$$

Multiplying both sides of the above equation by the matrix function $A(x)$, we see that $v \in W_0^{1,2}(B_R)$ solves the equation

$$\operatorname{div}(A\nabla v) = \operatorname{div}(AF) \qquad \text{in } B_R.$$

By Lemma 3.1, $v \in W_0^{1,p}(B_R)$ with

$$\|v\|_{L^p(B_R; dx)} \le c\|AF\|_{L^p(B_R; dx)} \le c_1 \|F\|_{L^p(B_R; dx)}.$$

If we take $v = 0$ on $B_R^c$, then $v \in W^{1,p}(\mathbb{R}^n)$ and so by the Sobolev embedding theorem [16], Theorem 7.10, $v \in C(\mathbb{R}^n)$. This completes the proof of the lemma.  $\square$

Henceforth, we select and fix a ball $B_R \supset D$.

THEOREM 3.3.   *Assume the Markovian condition (1.7) is satisfied, $p > n$ and that the diffusion matrix $A$ satisfies the condition (3.1) of Lemma 3.1. Assume that $|\widetilde{\mathbf{b}}| \in L^p(D; dx)$ , and that $\mathbf{1}_D(|\mathbf{b}|^2 + |q|) \in \mathbf{K}_n$. Let $u$ be the unique weak solution of the Dirichlet boundary value problem (1.8). Then for $y \in \partial D$ that is regular for $(\frac{1}{2}\Delta, D)$, we have*

$$\lim_{x \to y, x \in D} u(x) = f(y).$$

Proof. It is enough to prove the theorem for nonnegative function $f$. Let $\{\Omega, \mathcal{F}, X_t, \mathbf{Q}_x, x \in \mathbb{R}^n\}$ denote the diffusion process defined as in (2.8). As in the proof of Theorem 2.3, put

$$\widehat{M}_t = \int_0^{t \wedge \tau_{B_R}} (A^{-1}\widehat{\mathbf{b}})(X_s)\, dM_s.$$

By Lemma 3.2, there exits a bounded, function $v \in W_0^{1,p}(B_R) \subset W_0^{1,2}(B_R)$ such that

$$\widehat{M}_t \circ r_t = -\widehat{M}_t + N_t^v \qquad \text{for } t < \tau_{B_R}.$$

Note that $\widehat{M}$ is a MAF of $X$ in the strict sense and $N^v$ is a CAF of $X$ in the strict sense of zero energy in view of [13], Theorem 1, since $v \in W_0^{1,p}(B_R)$ and so it is continuous on $\mathbb{R}^n$ if we extend $v$ to take value 0 on $B_R^c$. Therefore $\widehat{M}_t \circ r_t$ can be refined to be a CAF of $X$ in the strict sense. It follows that $u$ can be expressed as

$$u(x) = \mathbf{E}_x^{\mathbf{Q}}[f(X_{\tau_D}) \exp(A_{\tau_D})],$$

where

$$A_t = N_t^v + \int_0^t q(X_s)\, ds.$$

Under condition (1.7), the CAF $A_t$ is negative and decreasing in $t \in [0, \tau_{B_R})$. Hence for $x \in D$,

$$|u(x) - \mathbf{E}_x^{\mathbf{Q}}[f(X_{\tau_D})]| \leq \|f\|_\infty \mathbf{E}_x^{\mathbf{Q}}[|\exp(A_{\tau_D}) - 1|] \leq \|f\|_\infty \mathbf{E}_x^{\mathbf{Q}}[1 - \exp(A_{\tau_D})].$$

We know from [8] that for $y \in \partial D$ that is regular for $(\frac{1}{2}\Delta, D)$,

$$\lim_{x \to y, x \in D} \mathbf{E}_x^{\mathbf{Q}}[f(X_{\tau_D})] = f(y)$$

and $\lim_{x \to y, x \in D} \mathbf{Q}_x(\tau_D > t) = 0$ for every $t > 0$. Thus, it suffices to show that

$$\lim_{t \downarrow 0} \lim_{x \to y, x \in D} \mathbf{E}_x^{\mathbf{Q}}[\exp(A_{\tau_D \wedge t})] = 1.$$

For this, note that by Jensen's inequality,

$$1 \geq \mathbf{E}_x^{\mathbf{Q}}[\exp(A_{\tau_D \wedge t})] \geq \exp(\mathbf{E}_x^{\mathbf{Q}}[A_{\tau_D \wedge t}]).$$

On the other hand, by the Cauchy–Schwarz inequality,

$$(\mathbf{E}_x^{\mathbf{Q}}[A_{\tau_D \wedge t}])^2 \leq \mathbf{E}_x[H_{\tau_D \wedge t}^2]\mathbf{E}_x[A_{\tau_D \wedge t}^2],$$

where $H$ is the martingale given by (2.9). We know from [8] that $\sup_{x \in D} \sup_{s \in [0,1]} \mathbf{E}_x[H_s^2] < \infty$. Observe that

$$A_t = N_t^v + \int_0^t q(X_s)\, ds = v(X_t) - v(X_0) - \int_0^t \nabla v(X_s)\, dM_s + \int_0^t q(X_s)\, ds.$$

Since $v$ is bounded and continuous, it is known (see, e.g., [8]) that

$$\lim_{x\to y, x\in D} \mathbf{E}_x[(v(X_{\tau_D\wedge t}))^2] = v(x)^2 \quad\text{and}\quad \lim_{x\to y, x\in D} \mathbf{E}_x[v(X_{\tau_D\wedge t})] = v(x).$$

Thus,

$$\lim_{x\to y, x\in D} \mathbf{E}_x[A_{\tau_D\wedge t}^2]$$

$$\leq \lim_{x\to y, x\in D} c\mathbf{E}_x\Big[\big(v(X_{\tau_D\wedge t}) - v(X_0)\big)^2$$

$$+ \int_0^{\tau_D\wedge t} |\nabla v(X_s)|^2\,ds + \Big(\int_0^{\tau_D\wedge t} |q(X_s)|\,ds\Big)^2\Big]$$

$$\leq c\lim_{x\to y, x\in D} \mathbf{E}_x\Big[\int_0^{\tau_D\wedge t} |\nabla v(X_s)|^2\,ds$$

$$+ 2\int_0^{\tau_D\wedge t} |q(X_s)|\Big(\int_s^{\tau_D\wedge t} |q(X_r)|\,dr\Big)ds\Big]$$

$$\leq c\lim_{x\to y, x\in D} \mathbf{E}_x\Big[\int_0^{\tau_D\wedge t} |\nabla v(X_s)|^2\,ds$$

$$+ 2\Big(\sup_{z\in D} G_D|q|(z)|\Big)\int_0^{\tau_D\wedge t} |q(X_s)|\,ds\Big]$$

$$= 0.$$

In the second to the last inequality we used the Markov property of $X$, while in the last equality we used the fact that $|\nabla v|^2$ and $q$ are in the Kato class and the dominated convergence theorem. This proves that $\lim_{t\downarrow 0}\lim_{x\to y, x\in D} \mathbf{E}_x^{\mathbf{Q}} \times [A_{\tau_D\wedge t}] = 0$. Consequently $\lim_{t\downarrow 0}\lim_{x\to y, x\in D} \mathbf{E}_x^{\mathbf{Q}}[\exp(A_{\tau_D\wedge t})] = 1$. This proves the theorem. $\square$

**4. Without Markovian assumption (1.7).** In this section, we will drop the Markovian condition (1.7) and give a general result on the existence and uniqueness of the weak solutions of the Dirichlet boundary value problem (1.8).

Let $\mathbf{h} = (h_1(x), \ldots, h_n(x)): \mathbb{R}^n \to \mathbb{R}^n$ be a measurable function such that $\mathbf{h} \in L^p(\mathbb{R}^n \to \mathbb{R}^n)$ for some $p > n$. Let $\mu$ be a signed measure in $\mathbf{K}_n$. Consider

$$\mathcal{G} = \frac{1}{2}\sum_{i,j=1}^n \frac{\partial}{\partial x_i}\Big(a_{ij}(x)\frac{\partial}{\partial x_j}\Big) + \sum_{i=1}^n h_i(x)\frac{\partial}{\partial x_i} + \mu.$$

The quadratic from $(\mathcal{C}, \mathcal{D}(\mathcal{C}))$ associated with $\mathcal{G}$ is given by $\mathcal{D}(\mathcal{C}) = W^{1,2}(\mathbb{R}^n)$ and for $u, v \in W^{1,2}(\mathbb{R}^n)$,

$$\mathcal{C}(u, v) = (-\mathcal{G}u, v)$$

$$(4.1) \qquad = \frac{1}{2} \sum_{i,j=1}^{n} \int_{\mathbb{R}^n} a_{ij}(x) \frac{\partial u}{\partial x_i} \frac{\partial v}{\partial x_j} \, dx - \sum_{i=1}^{n} \int_{\mathbb{R}^n} h_i(x) \frac{\partial u}{\partial x_i} v(x) \, dx$$

$$- \int_{\mathbb{R}^n} u(x) v(x) \mu(dx).$$

We will regard $(\mathcal{C}, W^{1,2}(\mathbb{R}^n))$ as a lower-order perturbation of the symmetric Dirichlet form $(\mathcal{E}, W^{1,2}(\mathbb{R}^n))$ associated with the infinitesimal generator $\frac{1}{2}\nabla(A\nabla)$. Recall that $X$ is the symmetric diffusion process with infinitesimal generator $\frac{1}{2}\nabla(A\nabla)$, or equivalently, with $(\mathcal{E}, W^{1,2}(\mathbb{R}^n))$. Let $A_t^\mu$ be the CAF of $X$ whose Revuz measure is $\mu$. It is proved in [8] (see also [20]) that the semigroup $\{T_t, t \geq 0\}$ associated with $(\mathcal{C}, W^{1,2}(\mathbb{R}^n))$ is given by

$$(4.2) \qquad \begin{aligned} T_t g(x) &= \mathbf{E}_x \bigg[ g(X_t) \exp \bigg( \int_0^t (A^{-1}\mathbf{h})(X_s) \, dM_s \\ &\qquad\qquad - \frac{1}{2} \int_0^t \mathbf{h} A^{-1} \mathbf{h}^*(X_s) \, ds + A_t^\mu \bigg) \bigg] \\ &= \mathbf{E}_x^{\mathbf{P}^*} [g(X_t) e^{A_t^\mu}], \end{aligned}$$

where $\mathbf{E}_x^{\mathbf{P}^*}$ stands for the expectation with respect to the diffusion measure $\{\mathbf{P}_x^*, x \in \mathbb{R}^n\}$, defined by

$$(4.3) \qquad \frac{d\mathbf{P}_x^*}{d\mathbf{P}_x}\bigg|_{\mathcal{F}_t} = H_t,$$

where

$$(4.4) \qquad H_t = \exp \bigg( \int_0^t (A^{-1}\mathbf{h})(X_s) \, dM_s - \frac{1}{2} \int_0^t \mathbf{h} A^{-1} \mathbf{h}^*(X_s) \, ds \bigg).$$

Let $\nu$ be a positive Radon measure of finite energy with respect to the symmetric Dirichlet form $(\mathcal{E}, W^{1,2}(\mathbb{R}^n))$. Note that due to the $L^p$-integrability of $\mathbf{h}$ and the Kato class condition on $\mu$, there exists a constant $\alpha_0 \geq 1$ such that for every $\alpha > \alpha_0$,

$$\mathcal{C}_\alpha(u, v) := \mathcal{C}(u, v) + \alpha(u, v)$$

is a positive definite quadratic form satisfying

$$c_\alpha^{-1} \mathcal{C}_\alpha(u, u) \leq \mathcal{E}_1(u, u) \leq c_\alpha \mathcal{C}_\alpha(u, u) \qquad \text{for every } u \in W^{1,2}(\mathbb{R}^n).$$

By Lax–Milgram theorem, for any $\alpha > \alpha_0$, there is a unique function in $D(\mathcal{E})$, denoted by $U_\alpha \nu$, such that

$$(4.5) \qquad \mathcal{C}_\alpha(U_\alpha \nu, v) = \int_{\mathbb{R}^n} v(x) \nu(dx)$$

for all $v \in D(\mathcal{E})$. $U_\alpha \nu$ is called the $\alpha$-potential of $\nu$ associated with the quadratic form $(\mathcal{C}_\alpha, D(\mathcal{C}))$. We have the following useful representation for $U_\alpha \nu$.

LEMMA 4.1.   *Let $A^\nu$ be the positive CAF of $X$ associated with the smooth measure $\nu$ in the strict sense. Then, for sufficiently large $\alpha$,*

$$U_\alpha \nu(x) = \mathbf{E}_x \left[ \int_0^\infty R_t^{(\alpha)} \, dA_t^\nu \right] = \mathbf{E}_x^{\mathbf{P}^*} \left[ \int_0^\infty \exp(-\alpha t + A_t^\mu) \, dA_t^\nu \right],$$

*where*

$$(4.6) \quad R_t^{(\alpha)} = \exp\left( \int_0^t (A^{-1} \mathbf{h})(X_s) \, dM_s - \frac{1}{2} \int_0^t \mathbf{h} A^{-1} \mathbf{h}^*(X_s) \, ds - \alpha t + A_t^\mu \right).$$

PROOF.   By the same proof of [14], Lemma 2.2.5, we see that there exists an $\mathcal{E}$-nest consisting of an increasing sequence $\{F_k, k \geq 1\}$ of compact sets of $\mathbb{R}^n$ such that $\|U_\alpha(\mathbf{1}_{F_k}\nu)\|_\infty < \infty$ for each $k \geq 1$. By restricting to $F_k$ if necessary, we may assume that $g(x) := U_\alpha \nu(x)$ is bounded in the sequel. By (4.5), it follows that for every $v \in W^{1,2}(\mathbb{R}^n)$,

$$
\begin{aligned}
\mathcal{E}(g, v) &= \frac{1}{2} \sum_{i,j=1}^n \int_{\mathbb{R}^n} a_{ij}(x) \frac{\partial g}{\partial x_i} \frac{\partial v}{\partial x_j} \, dx \\
(4.7) \qquad &= \sum_{i=1}^n \int_{\mathbb{R}^n} h_i(x) \frac{\partial g}{\partial x_i} v(x) \, dx + \int_{\mathbb{R}^n} g(x) v(x) \mu(dx) \\
&\quad + \int_{\mathbb{R}^n} v(x) \nu(dx) - \alpha \int_{\mathbb{R}^n} g(x) v(x) \, dx.
\end{aligned}
$$

By [14], Theorem 5.4.2, (4.7) implies that $\mathbf{P}_x$-a.s.

$$
\begin{aligned}
g(X_t) = g(X_0) + \int_0^t \nabla g(X_s) \, dM_s - \sum_{i=1}^n \int_0^t h_i(X_s) \frac{\partial g}{\partial x_i}(X_s) \, ds \\
- \int_0^t g(X_s) \, dA_s^\mu + \alpha \int_0^t g(X_s) \, ds - A_t^\nu.
\end{aligned}
$$

Note that $R_t^{(\alpha)}$ satisfies

$$dR_t^{(\alpha)} = R_t^{(\alpha)} (A^{-1} \mathbf{h})(X_t) \, dM_t + R_t^{(\alpha)} \, dA_t^\mu - \alpha R_t^{(\alpha)} \, dt.$$

By Itô's formula, it follows that

$$(4.8) \qquad g(X_t) R_t^{(\alpha)} = g(X_0) + K_t - \int_0^t R_s^{(\alpha)} A_s^\nu,$$

where

$$K_t := \int_0^t g(X_s) R_s^{(\alpha)} (A^{-1} \mathbf{h})(X_s) \, dM_s + \int_0^t R_s^{(\alpha)} \nabla g(X_s) \, dM_s$$

is a local martingale. Let $\{\tau_k, k \geq 1\}$ be an increasing sequence of stopping times approaching to $\infty$ such that $\{K_{t \wedge \tau_k}, t \geq 0\}$ is a martingale for every $k \geq 1$. It follows from (4.8) that

$$g(x) = \mathbf{E}_x \left[ \int_0^{t \wedge \tau_k} R_s^{(\alpha)} \, dA_s^\nu \right] + \mathbf{E}_x[g(X_{t \wedge \tau_k}) R_{t \wedge \tau_k}^{(\alpha)}].$$

Note that since $|\mu| \in \mathbf{K}_n$, it follows from Khasminskii's inequality, Cauchy–Schwarz inequality and [8], Theorem 3.1, $\{R_s^{(\alpha)}, 0 \leq s \leq T\}$ is uniformly integrable under $\mathbf{P}_x$ for every $x \in \mathbb{R}^n$ and $T > 0$. Letting $k \to \infty$, by dominated and monotone convergence theorems we get

$$(4.9) \qquad g(x) = \mathbf{E}_x \left[ \int_0^t R_s^{(\alpha)} \, dA_s^\nu \right] + \mathbf{E}_x[g(X_t) R_t^{(\alpha)}].$$

Since the Revuz measure $\mu$ of $A^\mu$ is in the Kato class, for sufficiently large $\alpha$, we have

$$\lim_{t \to \infty} \mathbf{E}_x[R_t^{(\alpha)}] = 0.$$

As $g$ is assumed to be bounded, letting $t \to \infty$ in (4.9) we obtain that

$$g(x) = \mathbf{E}_x \left[ \int_0^\infty R_s^{(\alpha)} A_s^\nu \right],$$

which completes the proof.   □

Assume from now on that $D$ is a bounded Lipschitz domain. In particular, the boundary of $D$ is regular with respect to $(\frac{1}{2}\Delta, D)$. Let $R$ be the multiplicative functional of $X$ defined as in (4.6) with $\alpha = 0$; that is,

$$R_t = \exp\left( \int_0^t (A^{-1}\mathbf{h})(X_s) \, dM_s - \frac{1}{2} \int_0^t \mathbf{h} A^{-1}\mathbf{h}^*(X_s) \, ds + A_t^\mu \right).$$

The following is a gauge theorem involving the Girsanov as well as Feynman–Kac functional, which extends the corresponding result in [15] where the diffusion matrix $A$ is the identity matrix and the Kato class measure $\mu$ is absolutely continuous with respect to the Lebesgue measure [i.e., $\mu(dx) = q(x) \, dx$].

THEOREM 4.2.   *Let $D$ be a bounded Lipschitz domain in $\mathbb{R}^n$, $\mathbf{h} \in L^p(\mathbb{R}^n \to \mathbb{R}^n)$ for some $p > n$ and $\mu$ a signed measure in $\mathbf{K}_n$. Define $u(x) := \mathbf{E}_x[R_{\tau_D}]$. Then if $u(x) < \infty$ for some $x \in D$, then $u$ is bounded between two positive constants in $D$.*

Proof. Let $G_D$ denote the Green function of the symmetric diffusion $X$ in $D$ and $\mu^+$, $\mu^-$ be the positive, negative part of $\mu$, respectively. Since $\mu$ is in the Kato class, we have by Jensen's inequality

$$\inf_{x \in D} u(x)$$

$$\geq \inf_{x \in D} \exp\left( \mathbf{E}_x \left[ \int_0^{\tau_D} (A^{-1}\mathbf{h})(X_s) \, dM_s - \frac{1}{2} \int_0^{\tau_D} \mathbf{h} A^{-1} \mathbf{h}^*(X_s) \, ds + A_{\tau_D}^\mu \right] \right)$$

$$= \inf_{x \in D} \exp\left( \mathbf{E}_x \left[ A_{\tau_D}^\mu - \frac{1}{2} \int_0^{\tau_D} \mathbf{h} A^{-1} \mathbf{h}^*(X_s) \, ds \right] \right)$$

$$\geq \exp(-\|G_D \mu^-\|_\infty - c\|G_D |h|^2\|_\infty).$$

Let $Y$ be the diffusion process with infinitesimal generator $\frac{1}{2} \sum_{i,j=1}^n \frac{\partial}{\partial x_i}(a_{ij}(x) \times \frac{\partial}{\partial x_j}) + \sum_{i=1}^n h_i(x) \frac{\partial}{\partial x_i}$. It is known (see [8]) that $Y$ can be obtained from $X$ through Girsanov transform by the exponential martingale $H$ in (4.4). That is, if for every $x \in \mathbb{R}^n$, we let $\mathbf{P}_x^*$ be the measure defined by

$$(4.10) \qquad \frac{d\mathbf{P}_x^*}{d\mathbf{P}_x}\bigg|_{\mathcal{F}_t} = H_t \qquad \text{for every } t \geq 0,$$

then the diffusion process $X$ under $\mathbf{P}_x^*$ has the same distribution as $Y$ starting from $x$.

Note that $R_{\tau_D} = H_{\tau_D} \exp(A_{\tau_D}^\mu)$. For every $k \geq 1$, we have by (62) of Sharpe [24],

$$\mathbf{E}_x^{\mathbf{P}^*}[k \wedge \exp(A_{\tau_D}^\mu)] = \mathbf{E}_x \left[ \int_0^{\tau_D} k \wedge \exp(A_t^\mu) \, d(-H_s) + (k \wedge \exp(A_{\tau_D}^\mu)) H_{\tau_D} \right]$$

$$= \mathbf{E}_x[(k \wedge \exp(A_{\tau_D}^\mu)) H_{\tau_D}].$$

Passing $k \to \infty$, we have

$$(4.11) \qquad \mathbf{E}_x^{\mathbf{P}^*}[\exp(A_{\tau_D}^\mu)] = \mathbf{E}_x[\exp(A_{\tau_D}^\mu) H_{\tau_D}] = \mathbf{E}_x[R_{\tau_D}].$$

Since $\mathbf{h} \in L^p(D; dx)$ for some $p > n$, by Ancona [1], the Green function of $Y$ in $D$ is comparable to that of $X$ in $D$. Hence $\mu$ is in the Kato class of $Y$ in $D$. By [5], Theorem 2.2, we have $x \mapsto \mathbf{E}_x^{\mathbf{P}^*}[\exp(A_{\tau_D}^\mu)]$ is either bounded or identically infinite on $D$. It then follows from (4.11) that $u(x)$ is either bounded on $D$ or identically infinite, which proves the theorem. $\quad\square$

Now consider the Dirichlet boundary value problem

$$(4.12) \qquad \mathcal{G}u = 0 \quad \text{in } D \quad \text{with} \quad u = f \quad \text{on } \partial D,$$

where $f \in C(\partial D)$. Recall that $\{\mathbf{P}_x^*, x \in D\}$ are the probability measures defined by (4.10).

THEOREM 4.3. *Recall that $R$ is the multiplicative functional of $X$ defined as in (4.6) with $\alpha = 0$. Assume that $\mathbf{E}_{x_0}[R_{\tau_D}] < \infty$ for some $x_0 \in D$. Then there exists a unique, continuous weak solution to the Dirichlet boundary value problem (4.12), which is given by*

$$(4.13) \qquad u(x) = \mathbf{E}_x[R_{\tau_D} f(X_{\tau_D})] = \mathbf{E}_x^{\mathbf{P}^*}[e^{A_{\tau_D}^\mu} f(X_{\tau_D})], \qquad x \in D.$$

PROOF. We first show that $u$ defined by (4.13) is a weak solution to the Dirichlet boundary value problem (4.12). Put

$$v_1(x) = \mathbf{E}_x^{\mathbf{P}^*}[f(X_{\tau_D})] \quad \text{and} \quad v_2(x) = \mathbf{E}_x^{\mathbf{P}^*}\left[\int_0^{\tau_D} v_1(X_s) e^{A_s^\mu} \, dA_s^\mu\right].$$

Then $u_1(x) = v_1(x) + v_2(x)$. By [8], Theorem 4.5, we know that $v_1$ is locally in $W^{1,2}(\mathbb{R}^n)$ and

$$\mathcal{C}(v_1, \phi) = -\int_D v_1(x) \phi(x) \mu(dx) \qquad \text{for every } \phi \in W_0^{1,2}(D).$$

Here the quadratic form $(\mathcal{C}, W^{1,2}(D))$ is defined by (4.1). We will show now that $v_2 \in W_0^{1,2}(D)$. Let

$$G_\beta g(x) = \mathbf{E}_x^{\mathbf{P}^*}\left[\int_0^{\tau_D} e^{-\beta t} e^{A_t^\mu} g(X_s) \, ds\right].$$

Denote by $\widehat{G}_\beta$ the adjoint operator of $G_\beta$ in $L^2(D; dx)$. Applying Lemma 4.1 to the process $X$ killed upon leaving $D$, we have

$$v_2(x) - \beta G_\beta^1 v_2(x) = \mathbf{E}_x^{\mathbf{P}^*}\left[\int_0^{\tau_D} e^{-\beta t} e^{A_t^\mu} v_1(X_t) \, dA_t^\mu\right] = \bar{U}_\beta(v_1\mu),$$

where $\bar{U}_\beta(v_1\mu)$ is defined as in the proof of Lemma 4.1 but with the killed process $X^D$ in place of $X$. Since $\mu$ belongs to the Kato class of $X$ (and of $Y$ by the proof of Theorem 4.2) in $D$, we have

$$\begin{aligned}
\beta(v_2 &- \beta G_\beta v_2, v_2) \\
&= \beta \int_D \mathbf{E}_x^{\mathbf{P}^*}\left[\int_0^{\tau_D} e^{-\beta t} e^{A_t^\mu} v_1(X_t) \, dA_t^\mu\right] v_2(x) \, dx \\
&= \beta \int_D \bar{U}_\beta(v_1\mu) v_2(x) \, dx \\
&= \int_D \beta \widehat{G}_\beta v_2(x) v_1(x) \mu(dx) \\
&\leq C_1 \left(\int_D |\nabla v_1|^2 \, dx + \int_D |v_1|^2 \, dx\right) \\
&\quad + C_2 \int_D |\beta \widehat{G}_\beta v_2|^2 \, dx + \mathcal{C}(\beta \widehat{G}_\beta v_2, \beta \widehat{G}_\beta v_2)
\end{aligned}$$

$$\leq C_1 \left( \int_D |\nabla v_1|^2 \, dx + \int_D |v_1|^2 \, dx \right) + C_2 \frac{\beta}{\beta - \alpha_0} \int_D |\beta \widehat{G}_\beta v_2|^2 \, dx$$
$$+ \frac{1}{2} \beta (v_2 - \beta G_\beta v_2, v_2),$$

where in the last inequality [19], Lemma 3.1(i), is used. Thus,

$$\sup_{\beta > 0} \beta (v_2 - \beta G_\beta v_2, v_2) < \infty,$$

which implies that $v_2 \in W_0^{1,2}(D)$. Moreover for $\phi \in W_0^{1,2}(D)$,

$$\mathcal{C}(v_2, \phi) = \lim_{\beta \to \infty} \beta(v_2 - \beta G_\beta v_2, \phi) = \lim_{\beta \to \infty} \int_D \beta \widehat{G}_\beta \phi(x) v_1(x) \mu(dx)$$
$$= \int_D v_1(x) \phi(x) \mu(dx).$$

The last equation follows from the fact that $\beta \widehat{G}_\beta \phi(x)$ converges to $\phi$ in the Dirichlet space $(\mathcal{E}_1, W_0^{1,2}(D))$. Thus, for $\phi \in W_0^{1,2}(D)$,

$$\mathcal{C}(u, \phi) = \mathcal{C}(v_1, \phi) + \mathcal{C}(v_2, \phi) = 0.$$

This means that $u$ is a weak solution to $\mathcal{G}u = 0$ in $D$. It is a well-known fact in the theory of PDE (cf. [16]) that $u$ is continuous inside $D$. Next we show that the boundary condition is fulfilled, that is,

$$\lim_{x \to y, x \in D} u_1(x) = f(y) \tag{4.14}$$

for every $y \in \partial D$. It is proved in [8] that every point in $\partial D$ is a regular point of $D^c$ with respect to $Y$

$$\lim_{x \to y, x \in D} v_1(x) = \lim_{x \to y, x \in D} \mathbf{E}_x^{\mathbf{P}^*}[f(X_{\tau_D})] = f(y)$$

for every $y \in \partial D$. As $u_1 = v_1 + v_2$, it suffices to show that

$$\lim_{x \to y, x \in D} v_2(x) = 0 \qquad \text{for } y \in \partial D.$$

Note that

$$v_2(x) = \mathbf{E}_x^{\mathbf{P}^*}[f(X_{\tau_D})(e^{A_{\tau_D}^\mu} - 1)].$$

For any $t > 0$, write $v_2$ as

$$v_2(x) = \mathbf{E}_x^{\mathbf{P}^*}[f(X_{\tau_D})(e^{A_{\tau_D}^\mu} - 1); \tau_D \leq t] + \mathbf{E}_x^{\mathbf{P}^*}[f(X_{\tau_D})(e^{A_{\tau_D}^\mu} - 1); \tau_D > t].$$

Now

$$|\mathbf{E}_x^{\mathbf{P}^*}[f(X_{\tau_D})(e^{A_{\tau_D}^\mu} - 1); \tau_D > t]| \leq \|f\|_\infty \mathbf{E}_x^{\mathbf{P}^*}[(e^{A_{\tau_D}^\mu} + 1); \tau_D > t]$$
$$= \|f\|_\infty \mathbf{E}_x^{\mathbf{P}^*}[e^{A_t^\mu} \mathbf{E}_{X_t}^{\mathbf{P}^*}[e^{A_{\tau_D}^\mu} + 1]; \tau_D > t]$$
$$\leq \|f\|_\infty \mathbf{E}_x^{\mathbf{P}^*}[e^{A_t^{|\mu|}}(M + 1); \tau_D > t].$$

Here in the last inequality, we used the assumption that $\mathbf{E}_x[R_{\tau_D}] < \infty$ for some $x \in D$ and thus by Theorem 4.2, $M := \sup_{x \in D} \mathbf{E}_x^{\mathbf{P}^*}[\exp(A_{\tau_D}^\mu] = \sup_{x \in D} \mathbf{E}_x[R_{\tau_D}] < \infty$. Let $y \in \partial D$. Since

$$\lim_{x \to y, x \in D} \mathbf{P}_x^*(\tau_D > t) = 0 \quad \text{and} \quad \sup_{x \in D} \mathbf{E}_x^{\mathbf{P}^*}[e^{A_t^{|\mu|}}] < \infty,$$

we have

$$\lim_{x \to y, x \in D} \mathbf{E}_x^{\mathbf{P}^*}[f(X_{\tau_D})(e^{A_{\tau_D}^\mu} - 1) : \tau_D > t] = 0.$$

Thus for every $t > 0$,

$$\begin{aligned}
\lim_{x \to y, x \in D} |v_2(x)| &= \limsup_{x \to y, x \in D} |\mathbf{E}_x^{\mathbf{P}^*}[f(X_{\tau_D})(e^{A_{\tau_D}^\mu} - 1); \tau_D \le t]| \\
&\le \|f\|_\infty \limsup_{x \to y, x \in D} \mathbf{E}_x^{\mathbf{P}^*}[e^{A_t^{|\mu|}} - 1; \tau_D \le t] \\
&\le \|f\|_\infty \limsup_{x \to y, x \in D} (\mathbf{E}_x[(H_t)^2])^{1/2}(\mathbf{E}_x[(e^{A_t^{|\mu|}} - 1)^2])^{1/2}.
\end{aligned}$$

By [8], Theorem 3.1, we know that $\sup_{x \in D} \mathbf{E}_x[(H_t)^2] < \infty$, where $H$ is defined in (4.10). On the other hand, since $\mu$ is in the Kato class of $X$, it follows from the Khasminskii's inequality that

$$\limsup_{t \to 0} \sup_{x \in D} \mathbf{E}_x[(e^{A_t^{|\mu|}} - 1)^2] = 0.$$

Thus we conclude that $\lim_{x \to y, x \in D} v_2(x) = 0$. So (4.14) is established.

To prove the uniqueness, let $u_1$ be any, bounded continuous weak solution of the Dirichlet boundary value problem (4.12). Then for $\phi \in W_0^{1,2}(D)$,

$$\begin{aligned}
(4.15) \qquad \mathcal{E}(u_1, \phi) &= \mathcal{C}(u_1, \phi) + \sum_{i=1}^n \int_D h_i(x) \frac{\partial u_1(x)}{\partial x_i} \phi(x) \, dx \\
&\quad + \int_D u_1(x) \phi(x) \mu(dx) \\
&= \sum_{i=1}^n \int_D h_i(x) \frac{\partial u_1(x)}{\partial x_i} \phi(x) \, dx + \int_D u_1(x) \phi(x) \mu(dx).
\end{aligned}$$

By [14], Theorem 5.4.2, it follows from (4.15) that the following decomposition holds:

$$\begin{aligned}
u_1(X_t) - u_1(X_0) &= \int_0^t \nabla u_1(X_s) \, dM_s - \sum_{i=1}^n \int_0^t h_i(X_s) \frac{\partial u_1}{\partial x_i}(X_s) \, ds \\
&\quad - \int_0^t u_1(X_s) \, dA_s^\mu
\end{aligned}$$

for $t < \tau_D$. Recall that

$$R_t = \exp\left( \int_0^t (A^{-1}\mathbf{h})(X_s)\, dM_s - \frac{1}{2} \int_0^t \mathbf{h} A^{-1} \mathbf{h}^*(X_s)\, ds + A_t^\mu \right)$$

satisfies

$$dR_t = R_t (A^{-1}\mathbf{h})(X_t)\, dM_t + R_t\, dA_t^\mu.$$

Applying Itô's formula, we get that for $t < \tau_D$,

$$(4.16) \quad d(u_1(X_t)R_t) = R_t \nabla u_1(X_t)\, dM_t + u_1(X_t)R_t(A^{-1}\mathbf{h})(X_t)\, dM_t.$$

This shows that $\{u_1(X_{t \wedge \tau_D})R_{t \wedge \tau_D}, t \geq 0\}$ is a $\mathbf{P}_x$-local martingale for every $x \in D$. We claim that $\{R_{t \wedge \tau_D}, t \geq 0\}$ is uniformly integrable with respect to $\mathbf{P}_x$ for every $x \in D$. To see this, write

$$R_{t \wedge \tau_D} = R_{\tau_D} \mathbf{1}_{\{t \geq \tau_D\}} + R_t \mathbf{1}_{\{t < \tau_D\}}.$$

Obviously the first term $\{R_{\tau_D} \mathbf{1}_{\{t \geq \tau_D\}}\}$ is uniformly integrable. For the second term, by Jensen's inequality,

$$\mathbf{1}_{\{t < \tau_D\}} \mathbf{E}_x[R_{\tau_D} | \mathcal{F}_t]$$
$$= \mathbf{1}_{\{t < \tau_D\}} \mathbf{E}_x[R_t(R_{\tau_D} \circ \theta_t) | \mathcal{F}_t]$$
$$= R_t \mathbf{1}_{\{t < \tau_D\}} \mathbf{E}_{X_t}[R_{\tau_D}]$$
$$\geq R_t \mathbf{1}_{\{t < \tau_D\}} \exp\left( -\frac{1}{2} \sup_{x \in D} \mathbf{E}_x\left[ \int_0^{\tau_D} \mathbf{h} A^{-1} \mathbf{h}^*(X_s)\, ds \right] + \inf_{x \in D} \mathbf{E}_x[A_{\tau_D}^\mu] \right)$$
$$\geq c R_t \mathbf{1}_{\{t < \tau_D\}}$$

for some positive constant $c$, where we have used the fact that $|h|^2, \mu$ are both in the Kato class. The above inequality implies that $\{R_t \mathbf{1}_{\{t < \tau_D\}}, t \geq 0\}$ is also $\mathbf{P}_x$-uniformly integrable. Therefore, $\{R_{t \wedge \tau_D}, t \geq 0\}$ is $\mathbf{P}_x$-uniformly integrable. Since $u$ is bounded continuous, we conclude that $\{u_1(X_{t \wedge \tau_D})R_{t \wedge \tau_D}, t \geq 0\}$ is uniformly integrable, hence a $\mathbf{P}_x$-martingale for every $x \in D$. Consequently,

$$\mathbf{E}_x[u_1(X_{t \wedge \tau_D})R_{t \wedge \tau_D}] = u_1(x).$$

Letting $t \to \infty$, we get that

$$u_1(x) = \mathbf{E}_x[f(X_{\tau_D})R_{\tau_D}],$$

which proves the uniqueness. $\square$

Now we can get to the main results of this section. Define

$$
\begin{aligned}
(4.17) \quad Z_t = \exp\Bigg\{ &\int_0^t (A^{-1}\mathbf{b})(X_s)\, dM_s + \left( \int_0^t (A^{-1}\widehat{\mathbf{b}})(X_s)\, dM_s \right) \circ r_t \\
&- \frac{1}{2} \int_0^t (\mathbf{b} - \widehat{\mathbf{b}}) A^{-1} (\mathbf{b} - \widehat{\mathbf{b}})^*(X_s)\, ds + \int_0^t q(X_s)\, ds \Bigg\}.
\end{aligned}
$$

Recall that $X$ is the symmetric diffusion with infinitesimal generator $\frac{1}{2}\nabla(A\nabla)$ and $M$ is the martingale part of $X$ in (1.9). For domain $D \subset \mathbb{R}^n$, $\tau_D := \inf\{t \geq 0 : X_t \notin D\}$ is the first exit time from $D$ by diffusion $X$. The following theorem is a new type of gauge theorem, in comparison with those found in [2, 5, 9].

THEOREM 4.4. *Let $D$ be a bounded Lipschitz domain contained in some ball $B_R$, $A$ be an $n \times n$ symmetric positive definitive matrix satisfying the condition (3.1) of Lemma 3.1, $|\mathbf{b}| + |\widehat{\mathbf{b}}| \in L^p(D; dx)$ for some $p > n$, and $\mathbf{1}_D q \in \mathbf{K}_n$. Then $Z_{\tau_D}$ of (4.17) is well defined under $\mathbf{P}_x$ for every $x \in D$. If $\mathbf{E}_{x_0}[Z_{\tau_D}] < \infty$ for some $x_0 \in D$, then the function $x \mapsto \mathbf{E}_x[Z_{\tau_D}]$ is bounded between two positive constants on $D$.*

PROOF. As before, put

$$\widehat{M}_t = \int_0^t (A^{-1}\widehat{\mathbf{b}})(X_s)\, dM_s \qquad \text{for } t \geq 0.$$

Let $R > 0$ so that $D \subset B_R := B(0, R)$. By Lemma 3.2, there exits a bounded function $v \in W_0^{1,p}(B_R) \subset W_0^{1,2}(B_R)$ such that

$$\widehat{M}_t \circ r_t = -\widehat{M}_t + N_t^v$$

and that $v \in W_0^{1,2}(B_R)$ satisfies the following equation in the distributional sense:

(4.18) $$\operatorname{div}(A\nabla v) = -2\operatorname{div}(\widehat{\mathbf{b}}) \qquad \text{in } B_R.$$

Note that by Sobolev embedding theorem, $v \in C(\mathbb{R}^n)$ if we extend $v = 0$ on $D^c$. Thus both $\widehat{M}$ and $N^v$ are CAFs of $X$ in the strict sense (cf. [13], Theorem 1), and so is $t \mapsto \widehat{M}_t \circ r_t$. Moreover,

$$\left(\int_0^{\tau_D} (A^{-1}\widehat{\mathbf{b}})(X_s)\, dM_s\right) \circ r_{\tau_D}$$

$$= -\int_0^{\tau_D} (A^{-1}\widehat{\mathbf{b}})(X_s)\, dM_s + N_{\tau_D}^v$$

$$= -\int_0^{\tau_D} (A^{-1}\widehat{\mathbf{b}})(X_s)\, dM_s + v(X_{\tau_D}) - v(X_0) - M_{\tau_D}^v$$

$$= -\int_0^{\tau_D} (A^{-1}\widehat{\mathbf{b}})(X_s)\, dM_s + v(X_{\tau_D}) - v(X_0)$$

$$\quad - \int_0^{\tau_D} \nabla v(X_s)\, dM_s.$$

Thus

$$(4.19) \quad \begin{aligned} Z_{\tau_D} = \frac{e^{v(X_{\tau_D})}}{e^{v(X_0)}} \exp\bigg(&\int_0^{\tau_D} (A^{-1}(\mathbf{b} - \widehat{\mathbf{b}}) - \nabla v)(X_s)\, dM_s \\ &+ \int_0^{\tau_D} \Big(q - \frac{1}{2}(\mathbf{b} - \widehat{\mathbf{b}}) A^{-1} (\mathbf{b} - \widehat{\mathbf{b}})^*\Big)(X_s)\, ds\bigg). \end{aligned}$$

So $Z_{\tau_D}$ is well defined under $\mathbf{P}_x$ for every $x \in D$. Since $v$ is bounded, $\mathbf{b} - \widehat{\mathbf{b}} \in L^p(D; dx)$, $\nabla v \in L^p(D; dx)$ and $\mathbf{1}_D q \in \mathbf{K}_n$, the theorem follows from Theorem 4.2. $\square$

Recall that $\mathcal{L}$ is the second-order differential operator defined by (1.3).

THEOREM 4.5. *Let $D$ be a bounded Lipschitz domain contained in some ball $B_R$, $A$ be an $n \times n$ symmetric positive definite matrix satisfying the condition (3.1) of Lemma 3.1, $|\mathbf{b}| + |\widehat{\mathbf{b}}| \in L^p(D; dx)$ for some $p > n$, and $\mathbf{1}_D q \in \mathbf{K}_n$. Let $Z$ be defined in (4.17) and assume that $\mathbf{E}_x[Z_{\tau_D}] < \infty$ for some $x \in D$. Then for every $f \in C(\partial D)$, there exists a unique weak solution $u$ to $\mathcal{L}u = 0$ in $D$ that is continuous on $\overline{D}$ with $u = f$ on $\partial D$. Moreover, the solution $u$ admits the following representation:*

$$(4.20) \quad u(x) = \mathbf{E}_x[Z_{\tau_D} f(X_{\tau_D})] \qquad \text{for } x \in D.$$

PROOF. Let $u$ be defined by the right-hand side of (4.20). Recall that $v$ is the function in $W_0^{1,p}(B_R)$ that is continuous on $\mathbb{R}^n$ if we extend $v = 0$ off $D$ in the proof of Theorem 4.4. Define

$$\begin{aligned} \widehat{Z}_{\tau_D} := \exp\bigg(&\int_0^{\tau_D} (A^{-1}(\mathbf{b} - \widehat{\mathbf{b}} - A\nabla v))(X_s)\, dM_s \\ &- \frac{1}{2} \int_0^{\tau_D} (\mathbf{b} - \widehat{\mathbf{b}} - A\nabla v) A^{-1} (\mathbf{b} - \widehat{\mathbf{b}} - A\nabla v)^*(X_s)\, ds \\ &- \int_0^{\tau_D} \langle \mathbf{b} - \widehat{\mathbf{b}}, \nabla v \rangle(X_s)\, ds \\ &+ \frac{1}{2} \int_0^{\tau_D} (\nabla v) A (\nabla v)^*(X_s)\, ds + \int_0^{\tau_D} q(X_s)\, ds\bigg). \end{aligned}$$

We have by (4.19) that

$$Z_{\tau_D} = \frac{e^{v(X_{\tau_D})}}{e^{v(X_0)}} \widehat{Z}_{\tau_D}.$$

Since the function $v$ is bounded and continuous on $\overline{D}$, it follows that $E_{x_0}[Z_{\tau_D}] < +\infty$ for some $x_0 \in D$ if and only if $E_{x_0}[\widehat{Z}_{\tau_D}] < +\infty$. Define $g := e^v u$. Then

$$(4.21) \quad g(x) = \mathbf{E}_x[\widehat{Z}_{\tau_D}(e^v f)(X_{\tau_D})] \qquad \text{for } x \in D.$$

Let

$$\mathcal{L}_1 := \frac{1}{2}\sum_{i,j=1}^n \frac{\partial}{\partial x_i}\Big(a_{ij}(x)\frac{\partial}{\partial x_j}\Big)$$

$$+ \sum_{i=1}^n (b_i(x) - \widehat{b}_i(x) - (A\nabla v)_i(x))\frac{\partial}{\partial x_i}$$

$$- \langle \mathbf{b} - \widehat{\mathbf{b}}, \nabla v \rangle(x) + \frac{1}{2}(\nabla v)A(\nabla v)^*(x) + q(x).$$

By Theorem 4.3, $g$ is a weak solution to the Dirichlet boundary value problem

$$\mathcal{L}_1 g = 0 \qquad \text{in } D \quad \text{and} \quad g = fe^v \qquad \text{on } \partial D.$$

Moreover, $g$ is continuous on $\overline{D}$ with $g = fe^v$ on $\partial D$. Hence, $u = e^{-v}g$ is bounded and continuous on $\overline{D}$ with $u = f$ on $\partial D$. Note that for any $\psi \in W_0^{1,2}(D)$,

$$\begin{aligned}
\mathcal{Q}^*(g,\psi) &:= (-\mathcal{L}_1 g, \psi) \\
&= \frac{1}{2}\sum_{i,j=1}^n \int_D a_{ij}\frac{\partial g}{\partial x_i}\frac{\partial \psi}{\partial x_j}\,dx \\
&\quad - \sum_{i=1}^n \int_D (b_i - \widehat{b}_i - (A\nabla v)_i)\frac{\partial g}{\partial x_i}\psi\,dx \\
&\quad - \int_D qg\psi\,dx + \int_{\mathbb{R}^n}\langle \mathbf{b} - \widehat{\mathbf{b}}, \nabla v\rangle g\psi\,dx \\
&\quad - \frac{1}{2}\int_D (\nabla v)A(\nabla v)^*(x)g\psi\,dx \\
&= 0.
\end{aligned} \tag{4.22}$$

Recalling the definition of the quadratic form $\mathcal{Q}$ from (1.6), we thus have for every $\phi \in C_c^\infty(D)$,

$$\begin{aligned}
\mathcal{Q}(u,\phi) &= \frac{1}{2}\sum_{i,j=1}^n \int_D a_{ij}\frac{\partial (ge^{-v})}{\partial x_i}\frac{\partial \phi}{\partial x_j}\,dx - \sum_{i=1}^n \int_D b_i\frac{\partial (ge^{-v})}{\partial x_i}\phi\,dx \\
&\quad - \sum_{i=1}^n \int_D \widehat{b}_i\frac{\partial \phi}{\partial x_i}ge^{-v}\,dx - \int_D qge^{-v}\phi\,dx \\
&= \frac{1}{2}\sum_{i,j=1}^n \int_D a_{ij}\frac{\partial g}{\partial x_i}\frac{\partial (e^{-v}\phi)}{\partial x_j}\,dx - \frac{1}{2}\sum_{i,j=1}^n \int_D a_{ij}\frac{\partial g}{\partial x_i}\frac{\partial e^{-v}}{\partial x_j}\phi\,dx
\end{aligned} \tag{4.23}$$

$$+ \frac{1}{2} \sum_{i,j=1}^{n} \int_{D} a_{ij} \frac{\partial e^{-v}}{\partial x_i} \frac{\partial \phi}{\partial x_j} g \, dx - \sum_{i=1}^{n} \int_{D} b_i \frac{\partial g}{\partial x_i} e^{-v} \phi \, dx$$

$$- \sum_{i=1}^{n} \int_{D} b_i \frac{\partial e^{-v}}{\partial x_i} g \phi \, dx - \sum_{i=1}^{n} \int_{D} \widehat{b}_i \frac{\partial \phi}{\partial x_i} g e^{-v} \, dx$$

$$- \int_{D} q g e^{-v} \phi \, dx.$$

Applying (4.22) with $\psi = \phi e^{-v}$ we obtain

$$\frac{1}{2} \sum_{i,j=1}^{n} \int_{\mathbb{R}^n} a_{ij}(x) \frac{\partial g}{\partial x_i} \frac{\partial (\phi e^{-v})}{\partial x_j} \, dx$$

$$= \sum_{i=1}^{n} \int_{\mathbb{R}^n} \left( b_i(x) - \widehat{b}_i(x) - (a \nabla v)_i(x) \right) \frac{\partial g}{\partial x_i} \phi e^{-v} \, dx$$

$$+ \int_{\mathbb{R}^n} q(x) g(x) \phi e^{-v} \, dx$$

$$- \int_{\mathbb{R}^n} \langle \mathbf{b} - \widehat{\mathbf{b}}, \nabla v \rangle(x) g(x) \phi e^{-v} \, dx$$

$$+ \frac{1}{2} \int_{\mathbb{R}^n} (\nabla v) A (\nabla v)^*(x) g \phi e^{-v} \, dx.$$

Substituting this expression into (4.23), we get after cancellations that

$$\mathcal{Q}(u, \phi) = - \sum_{i=1}^{n} \int_{\mathbb{R}^n} \widehat{b}_i(x) \frac{\partial (g\phi)}{\partial x_i} e^{-v} \, dx - \sum_{i=1}^{n} \int_{\mathbb{R}^n} (A \nabla v)_i(x) \frac{\partial g}{\partial x_i} \phi e^{-v} \, dx$$

$$+ \int_{\mathbb{R}^n} \langle \widehat{\mathbf{b}}, \nabla v \rangle \phi g e^{-v} \, dx + \frac{1}{2} \int_{\mathbb{R}^n} (\nabla v) A (\nabla v)^* g \phi e^{-v} \, dx$$

(4.24)

$$- \frac{1}{2} \sum_{i,j=1}^{n} \int_{\mathbb{R}^n} a_{ij}(x) \frac{\partial g}{\partial x_i} \frac{\partial e^{-v}}{\partial x_j} \phi \, dx$$

$$+ \frac{1}{2} \sum_{i,j=1}^{n} \int_{\mathbb{R}^n} a_{ij}(x) \frac{\partial e^{-v}}{\partial x_i} \frac{\partial \phi}{\partial x_j} g \, dx.$$

In the sequel, we write $\operatorname{div}(\cdot)$ for the divergence in the distribution sense. Now,

$$- \sum_{i=1}^{n} \int_{\mathbb{R}^n} \widehat{b}_i(x) \frac{\partial (g\phi)}{\partial x_i} e^{-v} \, dx$$

(4.25)

$$= \int_{\mathbb{R}^n} \operatorname{div}(\widehat{\mathbf{b}} e^{-v}) g \phi \, dx$$

$$= \int_{\mathbb{R}^n} \operatorname{div}(\widehat{\mathbf{b}}) e^{-v} g\phi \, dx - \int_{\mathbb{R}^n} \langle \widehat{\mathbf{b}}, \nabla v \rangle \phi g e^{-v} \, dx.$$

By virtue of (4.18),

$$\int_{\mathbb{R}^n} \operatorname{div}(\widehat{\mathbf{b}}) e^{-v} g\phi \, dx$$

$$= -\frac{1}{2} \int_{\mathbb{R}^n} \operatorname{div}(A\nabla v) e^{-v} g\phi \, dx$$

$$(4.26) \qquad = \frac{1}{2} \int_{\mathbb{R}^n} \langle A\nabla v, \nabla(e^{-v}g\phi) \rangle \, dx - \frac{1}{2} \int_{\mathbb{R}^n} \langle A\nabla v, \nabla v \rangle e^{-v} g\phi \, dx$$

$$+ \frac{1}{2} \int_{\mathbb{R}^n} \langle A\nabla v, \nabla g \rangle e^{-v} \phi) \, dx$$

$$+ \frac{1}{2} \int_{\mathbb{R}^n} \langle A\nabla v, \nabla \phi \rangle e^{-v} g\phi \, dx.$$

Combining (4.24)–(4.26), we arrive at

$$\mathcal{Q}(u,\phi) = -\int_{\mathbb{R}^n} (A(x)\nabla v(x) \cdot \nabla g(x)) \phi(x) e^{-v(x)} \, dx$$

$$- \frac{1}{2} \sum_{i,j=1}^{n} \int_{\mathbb{R}^n} a_{ij}(x) \frac{\partial g}{\partial x_i} \frac{\partial e^{-v}}{\partial x_j} \phi \, dx + \frac{1}{2} \int_{\mathbb{R}^n} \langle A\nabla v, \nabla g \rangle e^{-v} \phi \, dx$$

$$= 0.$$

This shows that $u$ is a weak solution to $\mathcal{L}u = 0$ in $D$. Recall that we have showed earlier that $u$ is continuous on $D$ with $u = f$ on $\partial D$.

It remains to show the uniqueness. Suppose that $u$ is a weak solution to $\mathcal{L}u = 0$ in $D$ and that $u$ is bounded and continuous on $\overline{D}$ with $u = f$ on $\partial D$. Then, for any $\phi \in W_0^{1,2}(D)$, $\mathcal{Q}(u,\phi) = 0$. Let $v$ be the solution of (4.18). Put $g(x) = e^{v(x)}u(x)$. Running the above proof backward, we can show that $\mathcal{Q}^*(g,\psi) = 0$ for any $\psi \in W_0^{1,2}(D)$. Hence $g$ is a weak solution to $\mathcal{L}_1 g = 0$ in $D$ that is bounded and continuous on $\overline{D}$ with $g = f e^v$ on $\partial D$. It follows from Theorem 4.3 that $g$ is given by (4.21). Arguing similarly as at the beginning of the proof, we conclude that $u$ can be expressed as in (4.20). This establishes the uniqueness. $\square$

REMARK 4.6. For $\widehat{\mathbf{b}} \in L^p(\mathbb{R}^n; dx)$, by [26], Theorem 4.1, there is a unique bounded weak solution $v \in W_0^{1,2}(B_R)$ such that

$$\operatorname{div}(A\nabla v) = -2\operatorname{div}(\widehat{\mathbf{b}}) \qquad \text{in } B_R.$$

Condition (3.1) of Lemma 3.1 is used to guarantee that $\nabla v \in L^p(B_R; dx)$ for some $p > n$ [see (4.18)], which in turn by the Sobolev embedding theorem implies that $v \in C(R^n)$ if we extend $v = 0$ on $D^c$. Condition (3.1) can be dropped from Theorems 4.4 and 4.5 if we assume directly that $\nabla v \in L^p(B_R; dx)$ for some $p > n$.

**5. Generalized Feynman–Kac transform.** In this section, we consider the special case of (1.3) with $\mathbf{b} = \hat{\mathbf{b}} = -A\nabla\rho$, for some $\rho \in W^{1,2}(\mathbb{R}^n)$ with $|\nabla\rho|^2 \in L^p(\mathbb{R}^n; dx)$ for some $p > n$. Note that by Sobolev embedding theorem, $\rho \in C(\mathbb{R}^n)$. As mentioned in the Introduction, the quadratic form $(\mathcal{Q}, W^{1,2}(\mathbb{R}^n))$ in (1.6) takes the following form:

$$
\begin{aligned}
(5.1) \quad \mathcal{Q}(u, v) = {} & \frac{1}{2} \sum_{i,j=1}^n \int_{\mathbb{R}^n} a_{ij}(x) \frac{\partial u}{\partial x_i} \frac{\partial v}{\partial x_j}\, dx + \frac{1}{2} \sum_{i,j=1}^n \int_{\mathbb{R}^n} a_{ij}(x) \frac{\partial(uv)}{\partial x_i} \frac{\partial \rho}{\partial x_j}\, dx \\
& + \int_{\mathbb{R}^n} u(x) v(x) q(x)\, dx
\end{aligned}
$$

for $u, v \in W^{1,2}(\mathbb{R}^n)$. Let $X$ be the symmetric diffusion with infinitesimal generator $\frac{1}{2}\nabla(A\nabla)$. When $q = 0$, it is established in Chen and Zhang [7] that $(\mathcal{Q}, W^{1,2}(\mathbb{R}^n))$ is the quadratic form associated with the generalized Feynman–Kac semigroup $\{T_t, t \geq 0\}$ defined by

$$
T_t f(x) = \mathbf{E}_x[e^{N_t^\rho} f(X_t)], \qquad x \in \mathbb{R}^n, t \geq 0,
$$

where $N^\rho$ is the CAF of $X$ in the Fukushima's decomposition of

$$
\rho(X_t) - \rho(X_0) = M_t^\rho + N_t^\rho, \qquad t \geq 0.
$$

(When $A$ is the identity matrix, the above result is proved in [17].) Recall that $M$ is the martingale defined in (1.9) and that $M_t^\rho = \int_0^t \nabla\rho(X_s)\, dM_s$. Note that (cf. [3])

$$
N_t^\rho = -\tfrac{1}{2}(M_t^\rho + M_t^\rho \circ r_t), \qquad t \geq 0.
$$

We thus have from (4.17) that, under the condition of Theorem 4.5, for every $f \in C(\partial D)$,

$$
u(x) := \mathbf{E}_x\left[\exp\left(N_{\tau_D}^\rho + \int_0^{\tau_D} q(X_s)\, ds\right) f(X_{\tau_D})\right], \qquad x \in \partial D,
$$

is the unique weak solution to $\mathcal{L}u = 0$ in $D$ that is continuous on $\overline{D}$ with $u = f$ on $\partial D$. However, we can do better in this case. More specifically, we can drop condition (3.1) of Lemma 3.1 in this case.

THEOREM 5.1. *Let $D$ be a bounded Lipschitz domain in $\mathbb{R}^n$. Let $A$ be a symmetric $n \times n$ bounded positive definite measurable matrix on $\mathbb{R}^n$, $\rho \in W^{1,2}(\mathbb{R}^n)$ with $\nabla\rho \in L^p(\mathbb{R}^n; dx)$ for some $p > n$. Let operator $\mathcal{L}$ be defined by (1.3) with $\mathbf{b} = \hat{\mathbf{b}} = -A\nabla\rho$ and $q \in \mathbf{K}_n$. Define*

$$
Z_t = \exp\left(N_t^\rho + \int_0^t q(X_s)\, ds\right), \qquad t \geq 0.
$$

*Then $Z_t$ is well defined under $\mathbf{P}_x$ for every $x \in \mathbb{R}^n$.*

(i) *If $\mathbf{E}_{x_0}[Z_{\tau_D}] < \infty$ for some $x_0 \in D$, then $x \mapsto \mathbf{E}_x[Z_{\tau_D}]$ is bounded between two positive constants on $D$.*

(ii) *Suppose that $\mathbf{E}_{x_0}[Z_{\tau_D}] < \infty$ for some $x_0 \in D$ and $f \in C(\partial D)$. Then*

$$u(x) := \mathbf{E}_x[Z_{\tau_D} f(X_{\tau_D})], \qquad x \in D,$$

*is the unique weak solution of $\mathcal{L}u = 0$ in $D$ that is continuous on $\overline{D}$ with $u = f$ on $\partial D$.*

PROOF. (i) Since $\nabla \rho \in L^p(\mathbb{R}^n; dx)$ for some $p > n$, by Sobolev embedding theorem, $\rho$ can be taken to be continuous on $\mathbb{R}^n$. By [13], Theorem 1, the following Fukushima decomposition in the strict sense holds:

$$\rho(X_t) - \rho(X_0) = M_t^\rho + N_t^\rho, \qquad t \geq 0,$$

where $M^\rho$ is the MAF of $X$ in the strict sense of finite energy and $N^\rho$ is the CAF of $X$ in the strict sense of zero energy. It follows

$$\begin{aligned}
Z_{\tau_D} &= \exp\left( \rho(X_{\tau_D}) - \rho(X_0) - M_{\tau_D}^\rho + \int_0^{\tau_D} q(X_s) \, ds \right) \\
&= \frac{\exp(\rho(X_{\tau_D}))}{\exp(u(X_0))} R_{\tau_D}
\end{aligned}$$

is well defined under $\mathbf{P}_x$ for every $x \in \mathbb{R}^n$, where

$$\begin{aligned}
(5.2) \quad R_t := \exp\Big( &- \int_0^t \nabla \rho(X_s) \, dM_s - \frac{1}{2} \int_0^t ((\nabla \rho) A (\nabla \rho)^*)(X_s) \, ds \\
&+ \int_0^t \left( \frac{1}{2} (\nabla \rho) A (\nabla \rho)^* + q \right)(X_s) \, ds \Big).
\end{aligned}$$

Since $\rho$ is bounded on $\overline{D}$ and $\nabla \rho \in L^p(\mathbb{R}^n; dx)$, $\mathbf{E}_{x_0}[Z_{\tau_D}] < \infty$ if and only if $\mathbf{E}_{x_0}[R_{\tau_D}] < \infty$, and so the conclusion of (i) follows from Theorem 4.2.

(ii) By (5.2),

$$e^{\rho(x)} u(x) = \mathbf{E}_x[R_{\tau_D}(e^\rho f)(X_{\tau_D})] \qquad \text{for } x \in D.$$

Define

$$\mathcal{G} := \frac{1}{2} \sum_{i,j=1}^n \frac{\partial}{\partial x_i}\left( (a_{ij}(x)) \frac{\partial}{\partial x_j} \right) - A(x)(\nabla \rho(x))^* \cdot \nabla + \frac{1}{2} (\nabla \rho) A (\nabla \rho)^* + q.$$

It follows from Theorem 4.3 that $v := e^\rho u$ is the unique weak solution for $\mathcal{G}v = 0$ in $D$ that is continuous on $\overline{D}$ with $v = e^\rho f$ on $\partial D$. Unwinding it for $u$ similar to that in the proof of Theorem 4.5, we conclude that $u \in W_{\text{loc}}^{1,2}(D)$ is the unique weak solution of $\mathcal{L}u = 0$ in $D$, that continuous on $\overline{D}$ with $u = f$ on $\partial D$. $\quad\square$

REMARK 5.2.    To better illustrate the main ideas of this paper, we have not strived to establish our theorems in the most general possible form. For example, in Theorems 4.4, 4.5 and 5.1, the potential function $q$ can be replaced by a signed Kato class measure $\mu$, just like in Theorems 4.2 and 4.3.

**Acknowledgments.**    We thank the anonymous referee and K. Kuwae for helpful comments on an earlier version of this paper.

DEPARTMENT OF MATHEMATICS                    SCHOOL OF MATHEMATICS
UNIVERSITY OF WASHINGTON                       UNIVERSITY OF MANCHESTER
SEATTLE, WASHINGTON 98195                      OXFORD ROAD
USA                                            MANCHESTER M13 9PL
E-MAIL: zchen@math.washington.edu              UNITED KINGDOM
                                               E-MAIL: tzhang@maths.manchester.ac.uk